\documentclass[11pt]{amsart}
\usepackage{amssymb,amsmath,bm,graphicx,float,cite,multirow,authblk}
\usepackage{subfigure,epstopdf,xcolor}
\usepackage{authblk}
\usepackage{graphicx}
\usepackage{lipsum}
\usepackage{algpseudocode}
\usepackage{multirow}
\usepackage{amsfonts}
\usepackage{subfigure}
\usepackage{epstopdf}
\usepackage{hyperref}
\usepackage{url}
\usepackage{color}
\usepackage{float}
\usepackage{tikz}
\usepackage{amsmath}  
\usepackage{algorithm}  
\usepackage{algpseudocode}
\usepackage{textcomp}
\usepackage{sidecap}
\usepackage{stmaryrd}

\usepackage{lipsum}

\theoremstyle{plain}

\newcommand{\bx}{\boldsymbol{x}}
\newcommand{\bu}{\boldsymbol{u}}
\newcommand{\bdf}{\boldsymbol{f}}
\newcommand{\bg}{\boldsymbol{g}}
\begin{document}

\author{Xurong Chi}
\address[X. Chi]{University of Science and Technology of China, Hefei, 230026, China}
\email{cxr123@mail.ustc.edu.cn}

\author{Jingrun Chen}
\address[J. Chen]{School of Mathematical Sciences, University of Science and Technology of China, Hefei, Anhui 230026, China; Suzhou Institute for Advanced Research, University of Science and Technology of China, Suzhou, Jiangsu 215123, China}
\email[Corresponding author]{jingrunchen@ustc.edu.cn}

\author{Zhouwang Yang}
\address[Z. Yang]{University of Science and Technology of China, Hefei, 230026, China}
\email{yangzw@ustc.edu.cn}

\title[RFM for Interface Problems]{The Random Feature Method for Solving Interface Problems}

\maketitle

\begin{abstract}
Interface problems have long been a major focus of scientific computing, leading to the development of various numerical methods. Traditional mesh-based methods often employ time-consuming body-fitted meshes with standard discretization schemes or unfitted meshes with tailored schemes to achieve controllable accuracy and convergence rate. Along another line, mesh-free methods bypass mesh generation but lack robustness in terms of convergence and accuracy due to the low regularity of solutions. In this study, we propose a novel method for solving interface problems within the framework of the random feature method. This approach utilizes random feature functions in conjunction with a partition of unity as approximation functions. It evaluates partial differential equations, boundary conditions, and interface conditions on collocation points in equal footing, and solves a linear least-squares system to obtain the approximate solution. To address the issue of low regularity, two sets of random feature functions are used to approximate the solution on each side of the interface, which are then coupled together via interface conditions. We validate our method through a series of increasingly complex numerical examples, including two-dimensional elliptic and three-dimensional Stokes interface problems, a three-dimensional elasticity interface problem, a moving interface problem with topological change, a dynamic interface problem with large deformation, and a linear fluid-solid interaction problem with complex geometry. Our findings show that despite the solution often being only continuous or even discontinuous, our method not only eliminates the need for mesh generation but also maintains high accuracy, akin to the spectral collocation method for smooth solutions. Remarkably, for the same accuracy requirement, our method requires two to three orders of magnitude fewer degrees of freedom than traditional methods, demonstrating its significant potential for solving interface problems with complex geometries.
\end{abstract}{\tiny }

\section{Introduction}
\label{sec1}
Interface problems are prevalent in various scientific and industrial applications, such as composite materials, multiphase flow and crystal growth.
These problems involve interfaces that divide an entire domain into several subdomains, which may have complex geometries and different physical properties.
Partial differential equations (PDEs) are often employed in the form of different (non)linear equations on different subdomains that are coupled together by interface conditions. As a result, solutions may exhibit non-smoothness or even discontinuity, making standard numerical methods for solving PDEs ineffective. In recent decades, significant progress has been made in studying interface problems from the algorithmic perspective. Among these studies, two kinds of methods have emerged, mesh-based methods and mesh-free methods.

Traditional mesh-based methods often employ time-consuming body-fitted meshes with standard discretization schemes or unfitted meshes with tailored schemes to achieve controllable accuracy and convergence rate. For body-fitted methods, meshes are generated to fit complex geometries without cutting through the interface. Adaptive refinement techniques \cite{chen2002efficiency, chen2004adaptive, chen2017interface} are  often employed, utilizing \textit{a posteriori} error indicators to guide adaptive mesh refinement and error control. However, constructing body-fitted meshes can be challenging and time-consuming for domains with complex geometries. In contrast, unfitted methods have significantly reduced the time required for mesh generation by allowing meshes to cut through the interface. Interface conditions are integrated into basis functions or tailored schemes. In the framework of finite difference method, there are the immersed interface method \cite{leveque1994immersed, li2006immersed}, the kernel-free boundary integral (KFBI) method \cite{DBLP:journals/corr/abs-2302-08022, DBLP:journals/corr/abs-2303-04992}, and the ghost fluid method \cite{liu2000boundary}. Meanwhile, various finite element-based methods have been utilized for interface problems. These include the partition of unity method, the generalized finite element method, the extended finite element method \cite{2003On, babuvska2017strongly, xiao2020high}, the penalty finite element method \cite{babuvska1970finite}, the matched interface and boundary method \cite{zhou2006high}, the hp-interface penalty finite element method \cite{wu2010unfitted}, and the unfitted discontinuous Galerkin methods \cite{massjung2012unfitted, bastian2009unfitted}. Among these methods, the immerse finite element method (IFM) deals with inhomogeneous interface conditions by constructing piecewise trilinear polynomials at each interface element \cite{li1998immersed, li2003new, gong2008immersed, chen2009adaptive}.

Along another line, mesh-free methods are proposed to avoid mesh generation. Radial basis function are commonly used for interface problems \cite{shin2011spectral, SIRAJULISLAM201838}, such as the generalized moving least squares approach \cite{dehghan2018interpolating,taleei2015efficient,zhang2008moving,trask2020compatible,hu2019spatially,trask2017high}. These methods incorporate techniques such as discontinuous derivative basis functions \cite{zhang2008moving}, gradient information \cite{trask2017high}, adaptive refinement strategies \cite{hu2019spatially}, and other techniques to improve the effectiveness. Recently, another popular type of mesh-free methods for solving interface problems is the machine learning-based method, leveraging their success in solving PDEs \cite{EHJ2017, HJE2018, EY2018, SS2018, PINN, Bao}. In \cite{wang2020mesh}, the non-smoothness of solution is addressed by employing piecewise neural networks, discretizing PDEs through collocation point sampling, and solving optimization problems using standard training algorithms like stochastic gradient descent. It shall be noted that while mesh-free methods eliminate the need for mesh generation, it is still common to truncate basis functions at the interface, which may encounter difficulty with complex geometries. In addition, the absence of convergence or convergence rate has been a serious issue for mesh-free methods.

The objective of this work is to propose a novel method for solving interface problems within the framework of the random feature method (RFM) \cite{JML-1-268}, which is a mesh-free method that has demonstrated its spectral accuracy on multiple complex problems. Our approach consists of three main components: we use random feature functions in conjunction with a partition of unity as approximation functions, evaluates PDEs, boundary conditions, and interface conditions on collocation points in equal footing, and solves a linear least-squares system to obtain the approximate solution. To address the issue of low regularity, two sets of random feature functions are used to approximate the solution on each side of the interface. This strategy allows us to tackle complex interface problems by utilizing collocation points instead of an underlying mesh, while achieves spectral accuracy in both stationary and time-dependent problems.

This paper is organized as follows. In Section \ref{sec2}, we define the interface problem and introduce the RFM framework, including the construction of approximate solutions, loss functions, and optimization procedures. Section \ref{sec3} presents numerical results of the proposed method for both stationary and time-dependent interface problems. In the first set of problems, the inhomogeneity of interface conditions is gradually increased. The remaining examples highlight the feasibility of RFM for solving time-dependent interface problems with complex geometries or intricate evolution. Comparison with available methods in the literature is included. Conclusions are drawn in Section \ref{sec4}.

\section{Model and Method}
\label{sec2}

This section will initially present the model of interface problem, followed by a sequential explanation of the approximation space, loss function, and optimization within the context of the RFM framework.

\subsection{Model of Interface Problem}
\label{sec2-1}
\begin{figure}[htbp]
	\centering
	\includegraphics[width=0.4\textwidth]{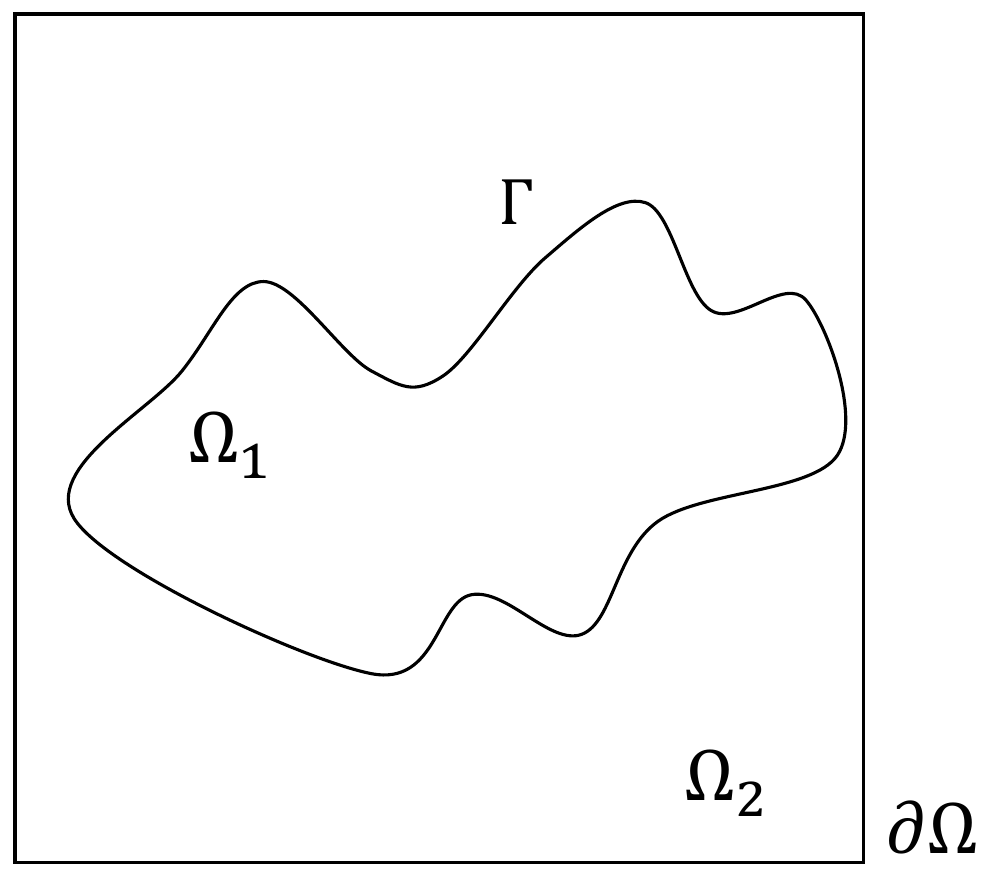}
	\caption{A sketch map for the domain $\Omega=\Omega_{1}\cup\Omega_{2}$ and the interface $\Gamma=\Omega_{1}\cap\Omega_{2}$.}
	\label{domain}
\end{figure}

Consider a bounded domain $\Omega \in \mathbb{R}^d$ that contains two subdomains $\Omega_{1}$ and $\Omega_{2}$ separated by a closed interface $\Gamma$, see Figure \ref{domain} for an illustration. 
It is a common assumption in traditional methods that the interface is $C^2$-continuous \cite{chen2021adaptive,2022High,xiao2020high,DBLP:journals/corr/abs-2302-08022}. 

In our formulation, we consider the following stationary interface problem \eqref{pde}
\begin{equation}
	\begin{cases}
		\mathcal{L}_{i}\bu(\boldsymbol{x}) =  \bdf_{i}(\boldsymbol{x}),  & \boldsymbol{x}\in\Omega_{i},\\
		\llbracket \bu(\boldsymbol{x}) \rrbracket=\boldsymbol{h}_1(\boldsymbol{x}), & \boldsymbol{x}\in \Gamma, \\
		\llbracket \boldsymbol{\sigma}(\boldsymbol{u}(\boldsymbol{x})) \boldsymbol{n} \rrbracket=\boldsymbol{h}_2(\boldsymbol{x}), & \boldsymbol{x}\in \Gamma,\\
		\mathcal{B}\bu(\boldsymbol{x}) =  \bg(\boldsymbol{x}),  & \boldsymbol{x}\in \partial\Omega,
	\end{cases}
	\label{pde}
\end{equation}
where operators $\mathcal{L}_{1}$ and $\mathcal{L}_{2}$ may have different forms or coefficients, $\bdf_{1}$, $\bdf_{2}$, $\boldsymbol{h}_1$, $\boldsymbol{h}_2$ and $\bg$ are given functions, $\boldsymbol{n}$ represents the outer unit normal vector of $\Gamma$ to the exterior of $\Omega_{1}$, and $\llbracket \boldsymbol{v} \rrbracket:=\boldsymbol{v}|_{\Omega_1}-\boldsymbol{v}|_{\Omega_2}$ denotes the jump of $\boldsymbol{v}$ across the interface $\Gamma$. $\boldsymbol{\sigma}$ is the stress tensor. When solving time-dependent equations, we treat time $t$ as an additional dimension as in \cite{EAJAM-13-435}.

\subsection{Approximation Space}
\label{sec2-2}
Following the random feature method \cite{JML-1-268}, we start by constructing the approximate solution $u_M$ as a linear combination of random feature functions
\begin{equation*}\label{representation1}
	u_M(\boldsymbol{x}) = \sum_{i=1}^M u_m \phi_m(\boldsymbol{x}).
\end{equation*}

\subsubsection{Random feature functions with partition of unity}
\label{sec2-2-1}
To capture the local variations near the interface, we construct multiple local solutions, each corresponding to a group of random feature functions, and piece them together using partition of unity (PoU) functions.

Specifically, we select $M_p$ points $\{\boldsymbol{x}_n\}_{n=1}^{M_p}$ from $\Omega$, and define the normalized coordinate
\begin{equation*}
	\boldsymbol{l}_{n}(\boldsymbol{x})=\frac{1}{\boldsymbol{r}_{n}}(\boldsymbol{x}-\boldsymbol{x}_{n}), \quad n=1,\cdots, M_p,
	\label{normalize}
\end{equation*}
which maps $\Omega^{(n)}=[x_{n1}-r_{n1},x_{n1}+r_{n1}]\times \cdots \times [x_{nd}-r_{nd},x_{nd}+r_{nd}]$ onto $[-1,1]^{d}$. Thus, a PoU function centered at $\boldsymbol{x}_n$ can be constructed. In one dimension, two commonly utilized types of PoU functions are given below.
\begin{tiny}
	\begin{equation*}
	\begin{aligned}
		& \psi^a_n(x)=\mathbb{I}_{[-1,1]}(l_n(x)), \\
		& \psi^b_n(x)=\mathbb{I}_{\left[-\frac{5}{4},-\frac{3}{4}\right]}(l_n(x)) \frac{1+\sin (2 \pi l_n(x))}{2}+\mathbb{I}_{\left[-\frac{3}{4}, \frac{3}{4}\right]}(l_n(x))+\mathbb{I}_{\left[\frac{3}{4}, \frac{5}{4}\right]}(l_n(x)) \frac{1-\sin (2 \pi l_n(x))}{2} .
		\label{pou}
	\end{aligned}
	\end{equation*}
\end{tiny}
Here $\psi^a_n(x)$ is discontinuous, while $\psi^b_n(x)$ is continuously differentiable. In high dimensions, the PoU function $\psi$ can be obtained directly from the tensor product $\psi_n(\boldsymbol{x})=\prod \limits_{k=1}^{d}\psi_n(x_{k})$. 

Then, $J_n$ random feature functions constructed by
\begin{equation*}
	\phi_{nj}(\boldsymbol{x}) = \sigma(\boldsymbol{k}_{nj} \cdot \boldsymbol{l}_{n}(\boldsymbol{x}) + b_{nj}), \quad j=1, \cdots, J_n,
\end{equation*}
where the nonlinear activation function $\sigma$ is often chosen as tanh or trigonometric functions and each component of $\boldsymbol{k}_{nj}$ and $b_{nj}$ is chosen uniformly from the interval $[-R_{nj},R_{nj}]$ and is fixed. In this way, the $n$-th locally space-dependent information is incorporated into $J_n$ random feature functions $\{\phi_{nj} \}_{j=1}^{J_n}$. Therefore, the degrees of freedom for the approximation space are given by $M= \sum_{n=1}^{M_p} J_n$.

A combination of these steps leads to the approximate solution $u_M$ as
\begin{equation}
	u_M(\boldsymbol{x})=\sum_{n=1}^{M_p} \psi_n (\boldsymbol{x})   \sum_{j=1}^{J_n }u_{nj} \phi_{nj} (\boldsymbol{x}).
	\label{representation2}
\end{equation}

\subsubsection{Approximate solution for interface problems}
\label{sec2-2-2}

To handle low regularity in interface problems, two sets of random features functions \eqref{representation2} are used to approximate the solution on each side of the interface, which are coupled together via interface conditions. This is implemented by constructing the approximate solution $u_M$ as
\begin{equation}
	u_M(\boldsymbol{x})=
	\begin{cases}
		\sum_{n=1}^{M_p} \psi_n (\boldsymbol{x}) \sum_{j=1}^{J_n }u_{nj}^{1} \phi_{nj}^{1} (\boldsymbol{x}),& \boldsymbol{x}\in\Omega_{1},\\
		\sum_{n=1}^{M_p} \psi_n (\boldsymbol{x}) \sum_{j=1}^{J_n }u_{nj}^{2} \phi_{nj}^{2} (\boldsymbol{x}),& \boldsymbol{x}\in\Omega_{2}.
	\end{cases}
	\label{representation3}
\end{equation}

For vectorial solutions, we approximate each component of the solution individually using \eqref{representation3}, i.e., 
\begin{equation*}\bu_M(\bx) = (u_M^{1}(\bx),\cdots,u_M^{K_I}(\bx))^T,\end{equation*}
where $K_I$ represents the output dimension.

By applying two sets of basis functions $\bu_i$ in $\Omega_{i}$, we reformaulate problem \eqref{pde} in the following form
\begin{equation}
	\begin{cases}
		\mathcal{L}_{1}\bu_{1}(\boldsymbol{x}) =  \bdf_{1}(\boldsymbol{x}),  & \boldsymbol{x}\in\Omega_{1},\\
		\mathcal{L}_{2}\bu_{2}(\boldsymbol{x}) =  \bdf_{2}(\boldsymbol{x}),  & \boldsymbol{x}\in\Omega_{2},\\
		\bu_{1}(\boldsymbol{x}) - \bu_{2}(\boldsymbol{x}) =\boldsymbol{h}_1(\boldsymbol{x}), & \boldsymbol{x}\in \Gamma, \\
		\boldsymbol{\sigma}(\bu_{1}(\boldsymbol{x}))\boldsymbol{n} - \boldsymbol{\sigma}(\bu_{2}(\boldsymbol{x})) \boldsymbol{n} =\boldsymbol{h}_2(\boldsymbol{x}), & \boldsymbol{x}\in \Gamma,\\
		\mathcal{B}\bu_{2}(\boldsymbol{x}) =  \bg(\boldsymbol{x}),  & \boldsymbol{x}\in \partial\Omega.
	\end{cases}
	\label{subdomain-pde}
\end{equation}

\subsection{Loss Function}
\label{sec2-3}
\begin{figure}[htbp]
	\centering
	\includegraphics[width=0.4\textwidth]{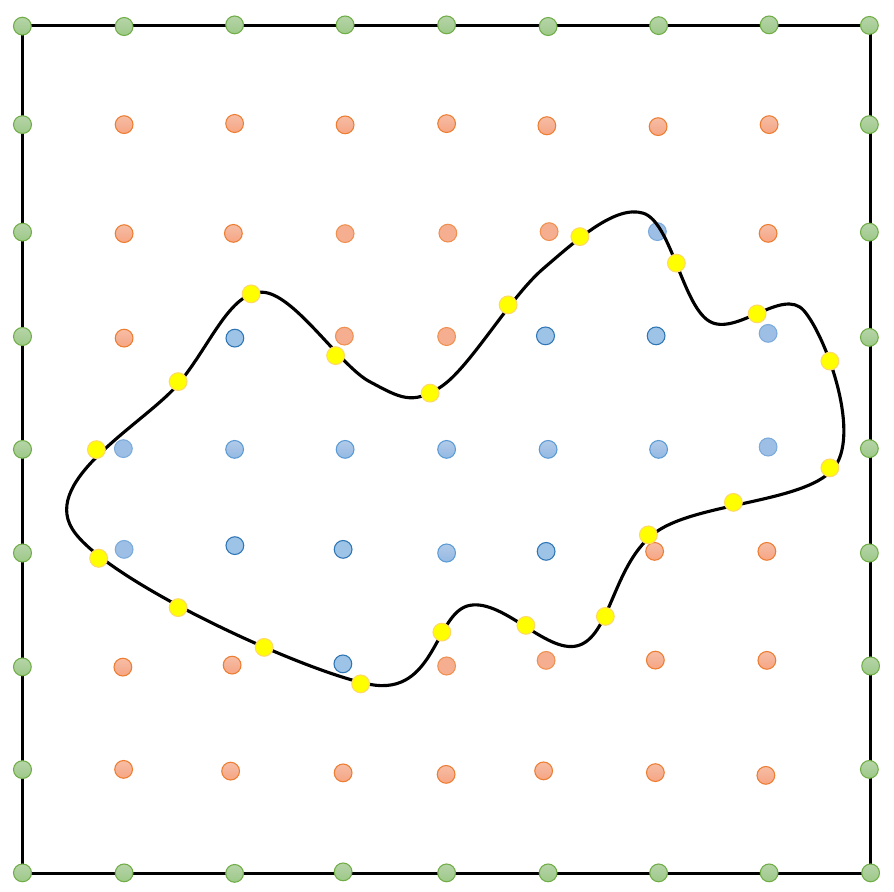}
	\caption{Collocation points for a two-dimensional domain: $C_{I}^{1}$, $\Omega_{1}$ interior points in orange; $C_{I}^{2}$, $\Omega_{2}$ interior points in blue; $C_{J}$, interface points in yellow; $C_{B}$, boundary points in green.}
	\label{collocation_points}
\end{figure}

Since the strong form is employed in the random feature method, we only requires the evaluation of PDE, interface condition, or boundary condition on collocation points to construct the loss function. Corresponding to problem \eqref{subdomain-pde}, we need to sample three sets of collocation points: $C_{I}=C_{I}^{1} \cup C_{I}^{2}$,  the set of interior points in $\Omega_{1} \cup \Omega_{2}$, $C_{J}$, the set of interface points on $\Gamma$ and $C_{B}$, the set of boundary points on $\partial\Omega$.
See Figure \ref{collocation_points} for an illustration.

Althrough there are many existing collocation point sampling methods, level set methods are widely used for interface tracing in interface problems \cite{saye2020review}. Therefore, we introduce the sign distance function $F(\boldsymbol{x})$ in level set methods to distinguish collocation points between two subdomains, as well as sample interface and boundary points.

Let $K_I$, $2K_J$ and $K_B$ be the number of conditions at each interior point, interface point and boundary point, respectively.
The total number of conditions is $N = K_I \# C_{I} +  2K_J \# C_{J} + K_B \# C_{B}$.
It is worth noting that, just like in random feature methods \cite{JML-1-268}, when $M_{p} > 1$ and PoU function $\psi_a$ is used, smoothness conditions between the adjacent elements in the partition are explicitly imposed by adding regularization terms in loss function \eqref{loss}, while no regularization is required when $\psi_b$ is used for second-order equations due to its first-order continuity.

Then, by introducing the penalty parameters $\boldsymbol{\lambda}_{Ii}^{1},\boldsymbol{\lambda}_{Ii}^{2} \in \mathbb{R}^{K_I \times 1}, \boldsymbol{\lambda}_{Jj}^{1},\boldsymbol{\lambda}_{Jj}^{2} \in \mathbb{R}^{K_J \times 1}, \boldsymbol{\lambda}_{Bk}\in \mathbb{R}^{K_B \times 1}$ for corresponding collocation points, we can construct the loss function in the following form
\begin{equation}
	\begin{aligned}
		Loss =& \sum_{\boldsymbol{x}_{i}^{1} \in C_{I}^{1}}\|\boldsymbol{\lambda}_{Ii}^{1}(\mathcal{L}_{1}\bu_{1}(\boldsymbol{x}_{i}^{1})-\bdf_{1}(\boldsymbol{x}_{i}^{1}))\|_{l^{2}}^{2}\\
		&+\sum_{\boldsymbol{x}_{i}^{2} \in C_{I}^{2}}\|\boldsymbol{\lambda}_{Ii}^{2}(\mathcal{L}_{2}\bu_{2}(\boldsymbol{x}_{i}^{2})-\bdf_{2}(\boldsymbol{x}_{i}^{2}))\|_{l^{2}}^{2}\\
		&+ \sum_{\boldsymbol{x}_{j} \in C_{J}}(\|\boldsymbol{\lambda}_{Jj}^{1}(\bu_{1}(\boldsymbol{x}_{j}) - \bu_{2}(\boldsymbol{x}_{j}) -\boldsymbol{h}_1(\boldsymbol{x}_{j}))\|_{l^{2}}^{2}\\
		&\quad\quad\quad + \|\boldsymbol{\lambda}_{Jj}^{2}(\boldsymbol{\sigma}(\bu_{1}(\boldsymbol{x}_{j}))\boldsymbol{n} - \boldsymbol{\sigma}(\bu_{2}(\boldsymbol{x}_{j})) \boldsymbol{n} -\boldsymbol{h}_2(\boldsymbol{x}_{j}))\|_{l^{2}}^{2})\\
		&+ \sum_{\boldsymbol{x}_{k} \in C_{B}}\|\boldsymbol{\lambda}_{Bk}(\mathcal{B}\bu_{2}(\boldsymbol{x}_{k})-\bg(\boldsymbol{x}_{k}))\|_{l^{2}}^{2}.
	\end{aligned}
	\label{loss}
\end{equation}

\subsection{Optimization}
\label{sec2-5}

Due to the significant differences in physical constants between two subdomains separated by the interface, it is important to balance the contributions from both subdomains in the loss function. 
The idea is to rescale each term in the loss function to the same order of magnitude based on the largest term in the sum.
Specifically, we choose the penalty parameters in loss function \eqref{loss} as follows:
\begin{equation*}
	\begin{aligned}
		&\lambda_{Ii}^{1\ell} = \frac{c}{\underset{1\leq n\leq M_p}{\max}\underset{1\leq j'\leq J_n}{\max}\underset{1\leq \ell'\leq K_I}{\max}|\mathcal{L}_{1}^{\ell} (\phi^{1\ell'}_{nj'}(\boldsymbol{x}_{i}^{1})\psi_{n}(\boldsymbol{x}_{i}^{1}))|},\quad \boldsymbol{x}_{i}^{1} \in C_{I}^{1}, \; \ell = 1,\cdots, K_I, \\
		&\lambda_{Ii}^{2\ell} = \frac{c}{\underset{1\leq n\leq M_p}{\max}\underset{1\leq j'\leq J_n}{\max}\underset{1\leq \ell'\leq K_I}{\max}|\mathcal{L}_{2}^{\ell} (\phi^{2\ell'}_{nj'}(\boldsymbol{x}_{i}^{2})\psi_{n}(\boldsymbol{x}_{i}^{2}))|},\quad \boldsymbol{x}_{i}^{2} \in C_{I}^{2}, \; \ell = 1,\cdots, K_I,\\
		&\lambda_{Jj}^{1\ell} = \frac{c}{\underset{m=1,2}{\max}\underset{1\leq n\leq M_p}{\max}\underset{1\leq j'\leq J_n}{\max}|\phi^{m\ell}_{nj'}(\boldsymbol{x}_{j})\psi_{n}(\boldsymbol{x}_{j})|},\quad \boldsymbol{x}_{j} \in C_{J}, \; \ell = 1,\cdots, K_J,\\
		&\lambda_{Jj}^{2\ell} = \frac{c}{\underset{m=1,2}{\max}\underset{1\leq n\leq M_p}{\max}\underset{1\leq j'\leq J_n}{\max}\underset{1\leq \ell'\leq K_I}{\max}|\boldsymbol{\sigma}^{\ell} (\phi^{m\ell'}_{nj'}(\boldsymbol{x}_{j})\psi_{n}(\boldsymbol{x}_{j}))\boldsymbol{n}(\boldsymbol{x}_{j})|},\quad \boldsymbol{x}_{j} \in C_{J}, \; \ell = 1,\cdots, K_J, \\
		&\lambda_{Bk}^{\ell} = \frac{c}{\underset{1\leq n\leq M_p}{\max}\underset{1\leq j'\leq J_n}{\max}\underset{1\leq \ell'\leq K_I}{\max}|\mathcal{B}^{\ell} (\phi^{\ell'}_{nj'}(\boldsymbol{x}_{k})\psi_{n}(\boldsymbol{x}_{k}))|},\quad \boldsymbol{x}_{k} \in C_{B}, \; \ell = 1,\cdots, K_B,
	\end{aligned}
\end{equation*}
where $c$ is a universal constant and we set $c=100$ in all experiments.

Since only the coefficients $\{u_{nj}^{1}\}\cup\{u_{nj}^{2}\}$ of the linear combination of basis functions are adjustable, the optimization problem is convex and can be solved using standard algorithms for linear least-squares approximation. Specifically, we use the \emph{SparseMatrix} class for matrix storage and the \emph{SPQR} solver to handle linear least-squares approximation in RFM. All implementations are based on the \emph{Eigen} \cite{eigenweb} which is a C++ template library for linear algebra.

\section{Numerical Results}
\label{sec3}

This section numerically examines two categories of interface problems. Initially, we employ stationary problems with progressively escalating inhomogeneity of interface conditions to numerically validate the convergence attributes of the RFM. Subsequently, we showcase results for time-dependent interface problems characterized by intricate evolution or complex geometry, thereby demonstrating the practicality of the RFM.

In all our experiments, unless otherwise specified, we use the default setup where the weights $\{k_m\}$ and $\{b_m\}$ are assumed to follow the distribution $\mathbb{U}[-1,1]$, the activation function is chosen as $\tanh$, the PoU function is $ \{\psi^{a}\}$.
We start by selecting a set of points $\{\boldsymbol{x}_{n}\}_{n=1}^{M_p}$ and construct the PoU functions. 
For each $\boldsymbol{x}_{n}$, we construct $J_{n}$ random feature functions with radius $\boldsymbol{r}_{n}$. 
Then we sample $Q$  equally spaced collocation points and use $F(\boldsymbol{x})$ to distinguish between two subdomains, where $\Omega_{1}=\{\boldsymbol{x}\in\Omega|F(\boldsymbol{x})\leq0\}$ and $\Omega_{2}=\{\boldsymbol{x}\in\Omega|F(\boldsymbol{x})\geq0\}$.
Boundary and interface collocation points are sampled using the level set method. This process results in a linear system $Au = b$, where $A$ is an $N \times M$ matrix. By default, all the results are measured in the relative $L^2$ error with a refined resolution. For time-dependent equations, the errors are evaluated at the final time.

\subsection{Stationary Problems}
\label{sec3-1}
We start with several stationary examples to investigate how the performance of RFM, specifically its accuracy, is affected by the inhomogeneity of interface conditions. Hyper-parameters in RFM for stationary problems are listed in Appendix \ref{sec::appendix}. 

\subsubsection{Elliptic interface problem with $h_1(\boldsymbol{x})=0$ and $h_2(\boldsymbol{x})=0$}
\label{sec3-1-1}
We first consider the elliptic interface problem as in \cite{chen2021adaptive,chen2022arbitrarily}
\begin{equation}
	\begin{cases}
		\begin{aligned}
			-\operatorname{div}(a(\boldsymbol{x}) \nabla u(\boldsymbol{x}))&=f(\boldsymbol{x}), & & \boldsymbol{x}\in \Omega_{1} \cup \Omega_{2}, \\
			\llbracket u(\boldsymbol{x}) \rrbracket&=0,& & \boldsymbol{x}\in \Gamma,\\
			\llbracket a(\boldsymbol{x}) \nabla u(\boldsymbol{x}) \cdot \boldsymbol{n}(\boldsymbol{x}) \rrbracket&=0, & & \boldsymbol{x}\in \Gamma, \\
			u(\boldsymbol{x})&=g(\boldsymbol{x}), & & \boldsymbol{x}\in  \partial \Omega.
		\end{aligned}
	\end{cases}
	\label{pde-elliptic}
\end{equation}

It is important to note that there are no jumps across the interface, i.e., $h_1(\boldsymbol{x})=h_2(\boldsymbol{x})=0$, but the diffusion coefficient can be highly contrasted.

We perform the calculation on the square domain $(-2.0,2.0)^2$ with a circular interface $F(\boldsymbol{x})=x^2+y^2-r^2=0$ of radius $r=1.1$. The exact solution and coefficient function are chosen to be
\begin{equation*}
	u(\boldsymbol{x})= \begin{cases}e^{|\boldsymbol{x}|^2-r^2}+10 r^2-1, & \boldsymbol{x} \in \Omega_{1} \\ 10|\boldsymbol{x}|^2, & \boldsymbol{x} \in \Omega_{2}\end{cases},
\end{equation*}
\begin{equation*}
	a(\boldsymbol{x})= \begin{cases}10, & \boldsymbol{x} \in \Omega_{1} \\ 1, & \boldsymbol{x} \in \Omega_{2}\end{cases}.
\end{equation*}

\begin{figure}[htbp]
	\centering
	\subfigure[$u$]{
		\includegraphics[width=0.28\textwidth]{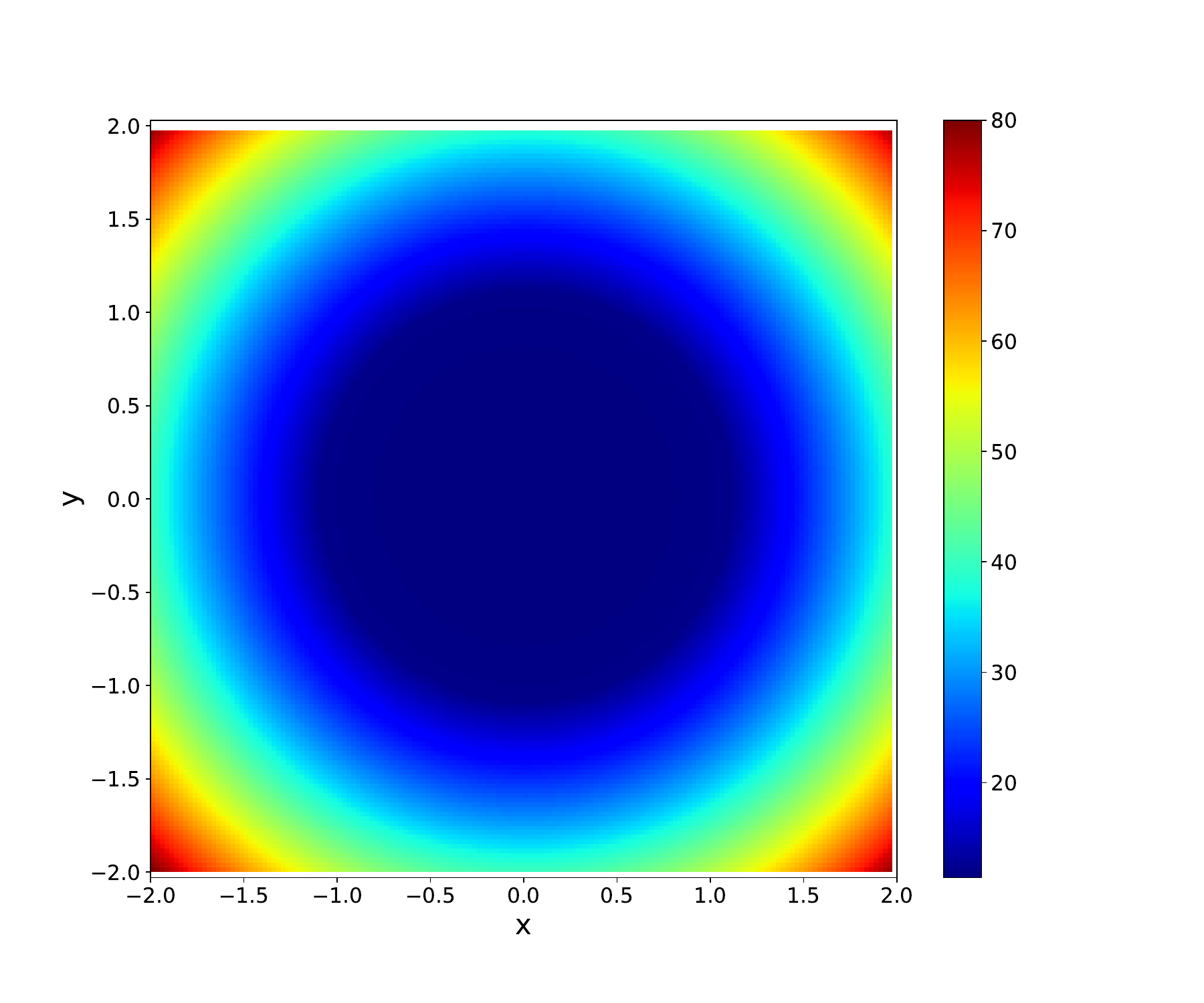}
	}
	\quad
	\subfigure[$u_x$]{
		\includegraphics[width=0.28\textwidth]{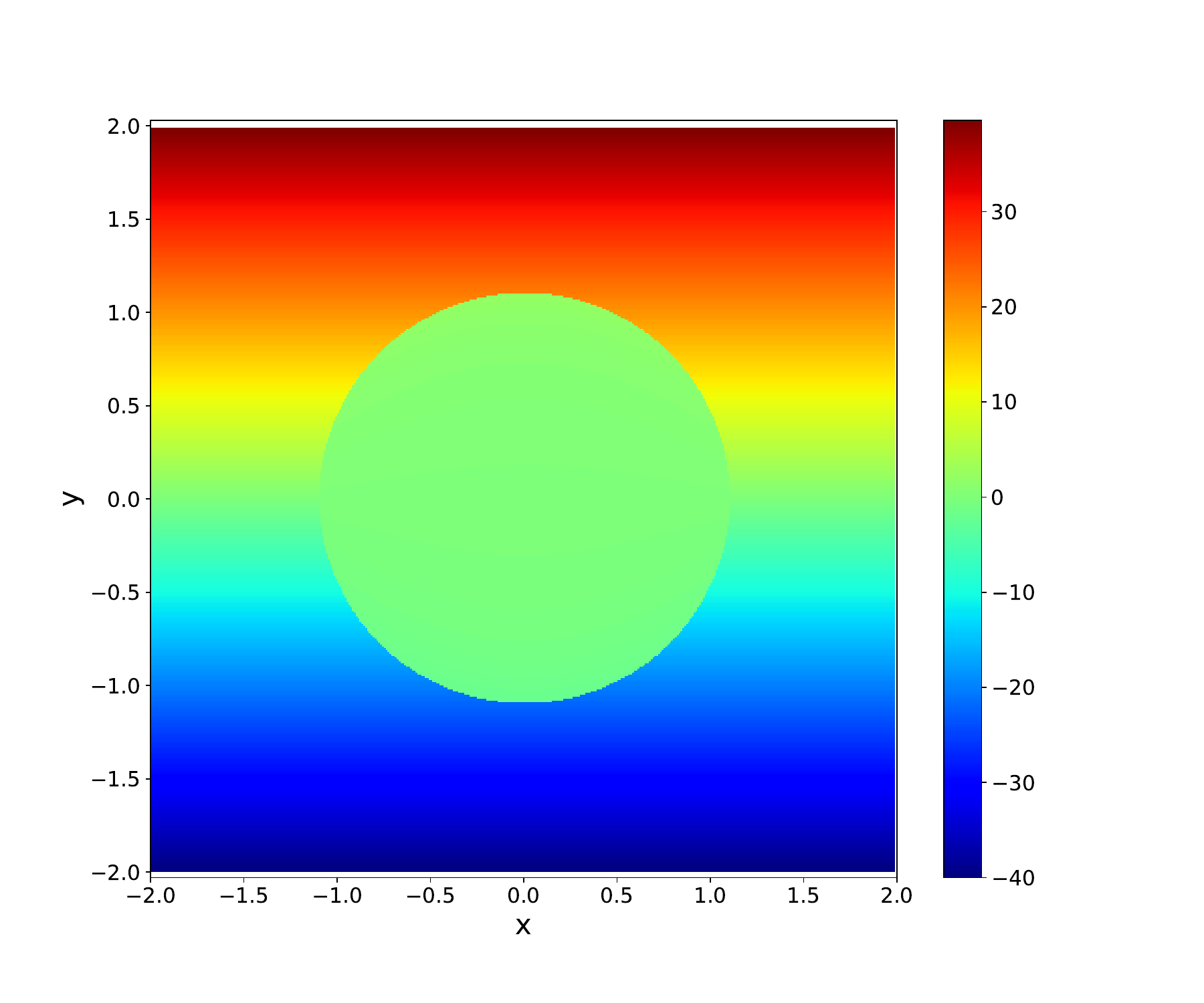}
	}
	\quad
	\subfigure[$u_y$]{
		\includegraphics[width=0.28\textwidth]{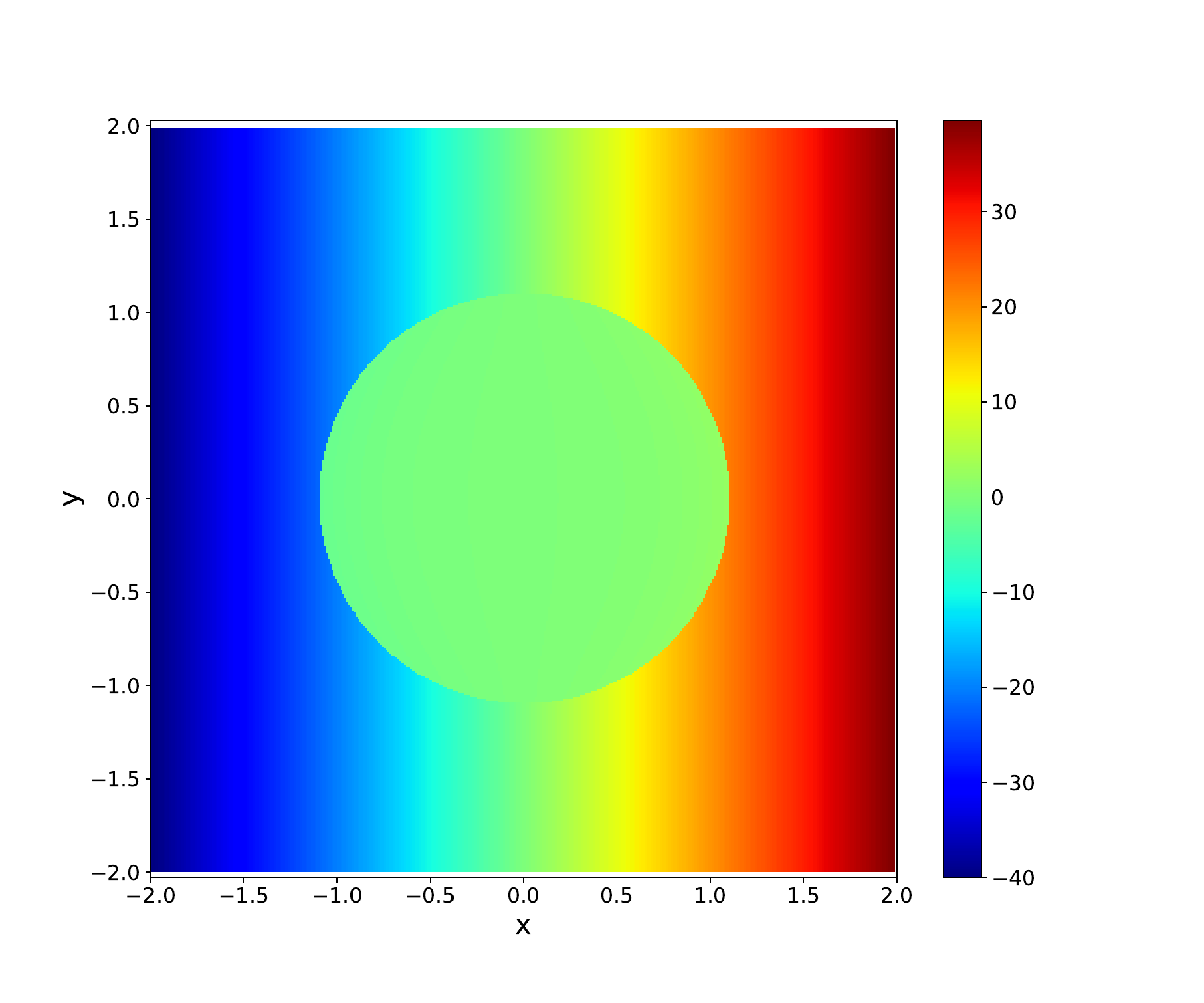}
	}
	\caption{Numerical solution and its first-order derivatives for two-dimensional elliptic interface problem \eqref{pde-elliptic}.}
	\label{figure3-1-1-1}
\end{figure}

Figure \ref{figure3-1-1-1} visualizes the numerical solution and its first-order derivatives obtained by RFM.

\begin{table}[htbp]
	\caption{\label{table3-1-1-1} RFM results for the two-dimensional elliptic interface equation.} \centering
	\begin{small}
		\begin{tabular}{|c|c|ccc|}
			\hline
			$M$ & $N$ & $u$ error & $u_{x}$ error & $u_{y}$ error\\
			\hline
			\multirow{3}{*}{1600} 
			& 3084  &  3.31E-5  &  2.96E-5  &  3.13E-5 \\
			& 9370  &  1.03E-8  &  2.28E-8  &  2.34E-8 \\
			& 18854 &  8.53E-10 &  4.71E-9  &  4.98E-9 \\
			\hline
			\multirow{3}{*}{3200} 
			& 3084  &  3.94E-6  &  3.36E-6  &  3.32E-6  \\
			& 9370  &  4.44E-9  &  3.81E-9  &  3.72E-9  \\
			& 18854 &  5.60E-10 &  6.33E-10 &  6.15E-10 \\
			\hline
		\end{tabular}
	\end{small}
\end{table}

Table \ref{table3-1-1-1} records the convergence of RFM, which exhibits the same exponential convergence observed in non-interface problems.
We attribute this to the fact that the exact solution is piecewisely smooth on each side of the interface and two sets of random feature functions are empolyed accordingly.

\begin{figure}[htbp]
	\centering
	\includegraphics[width=0.6\textwidth]{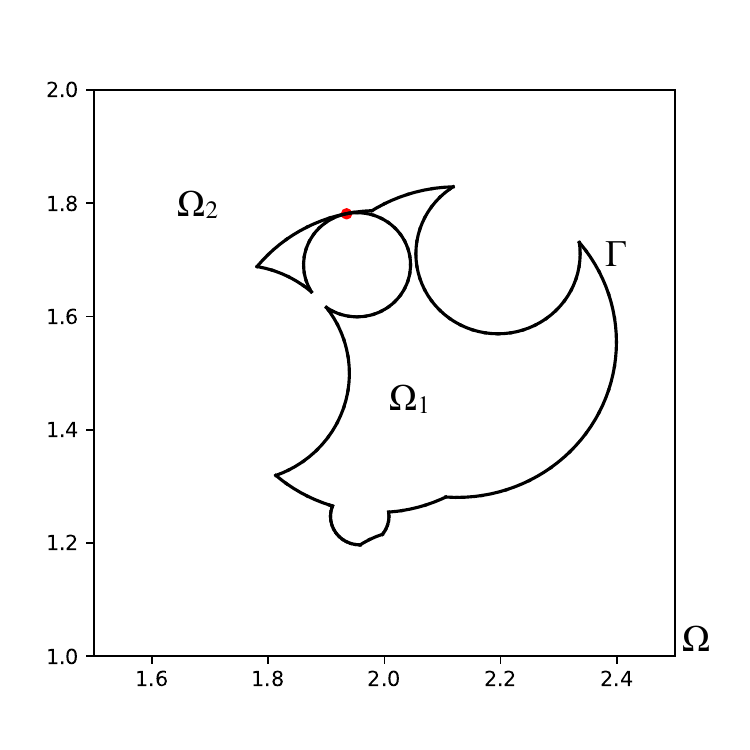}
	\caption{A two-dimensional complex interface.}
	\label{figure3-1-1-2}
\end{figure}

\begin{figure}[htbp]
	\centering
	\subfigure[$u$]{
		\includegraphics[width=0.28\textwidth]{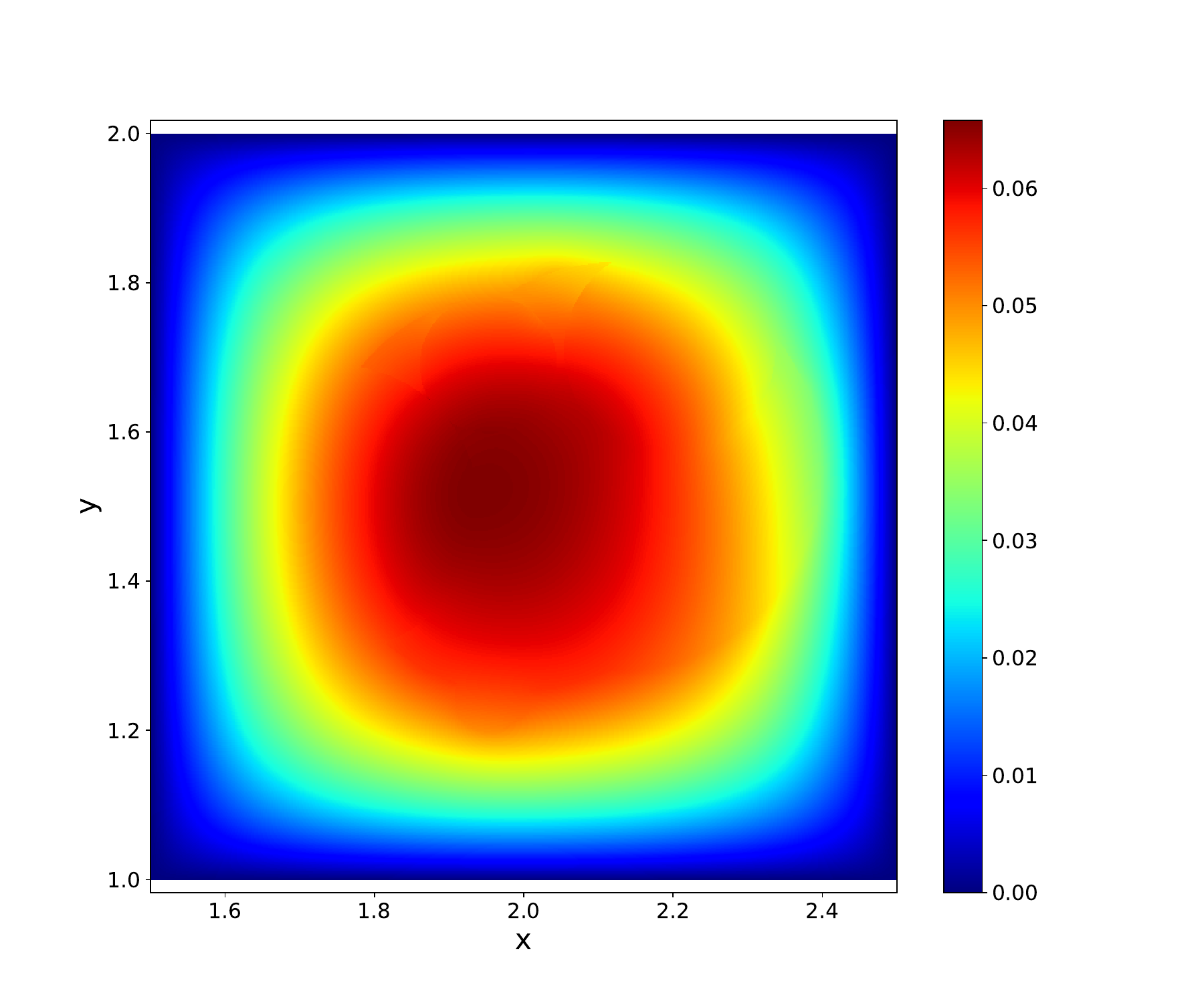}
	}
	\quad
	\subfigure[$u_x$]{
		\includegraphics[width=0.28\textwidth]{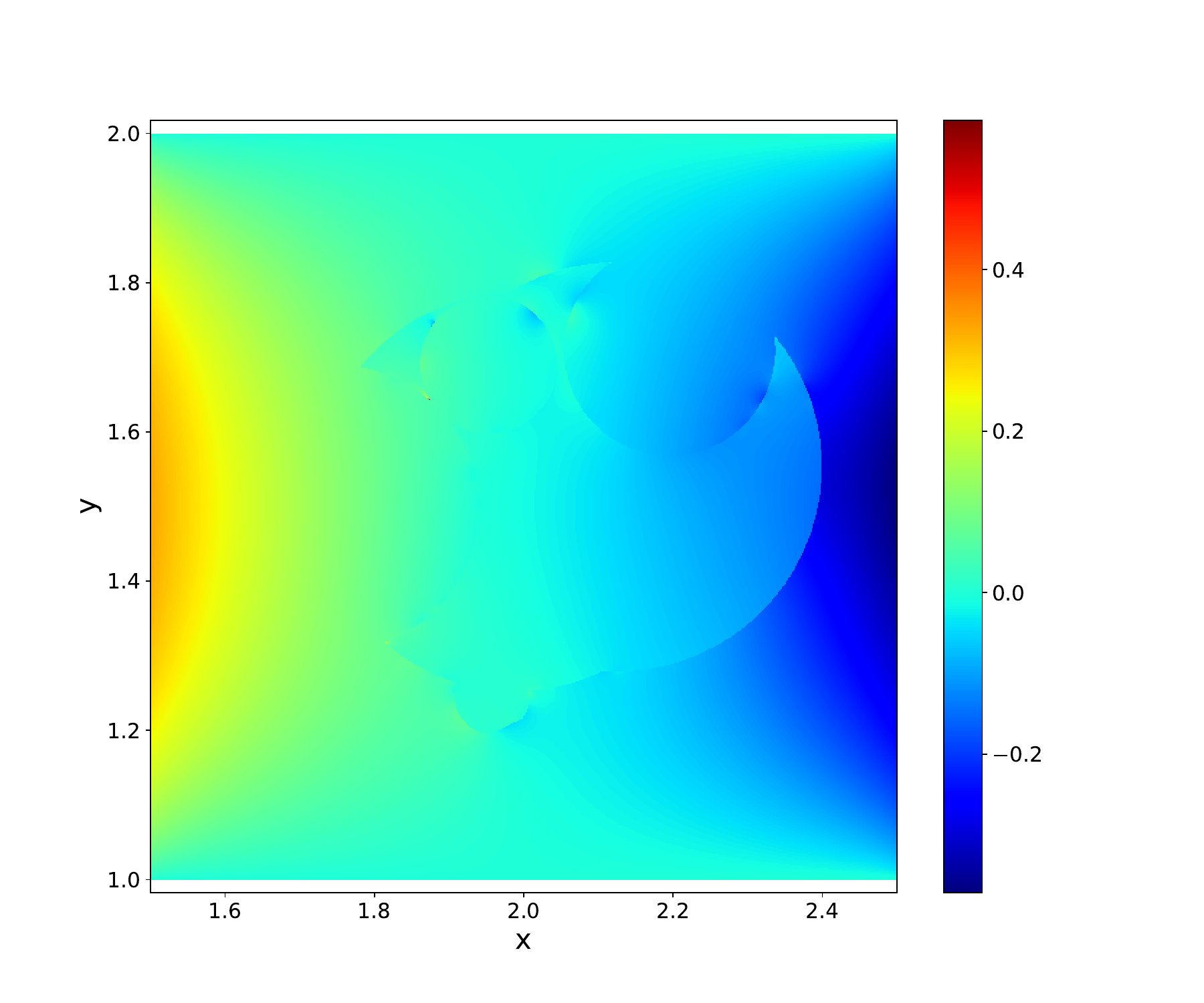}
	}
	\quad
	\subfigure[$u_y$]{
		\includegraphics[width=0.28\textwidth]{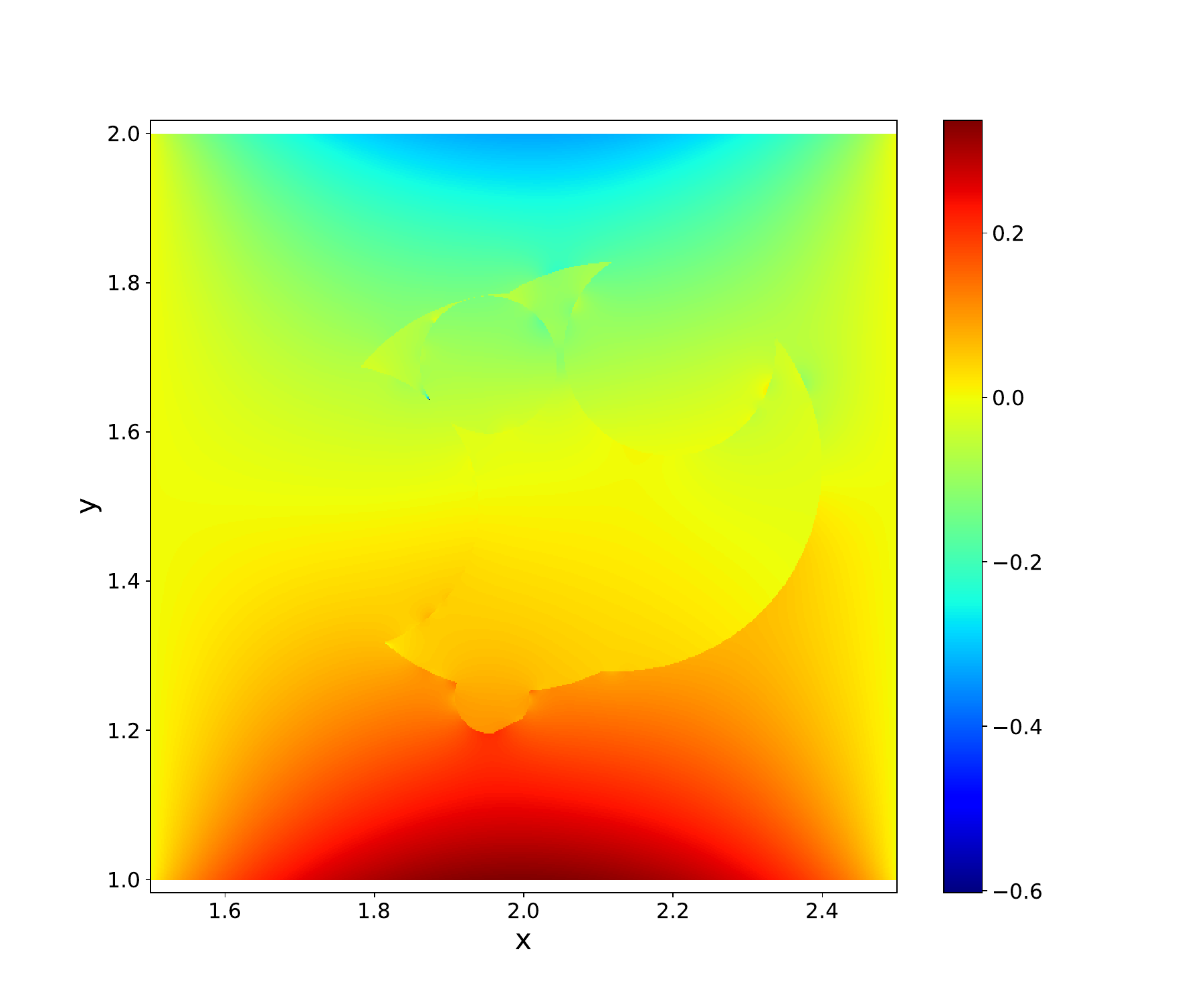}
	}
	\caption{Numerical solution and its first-order derivatives for two-dimensional elliptic interface problem with complex geometry in Figure \ref{figure3-1-1-2}.}
	\label{figure3-1-1-3}
\end{figure}

Further, we consider a complex interface geometry shown in Figure \ref{figure3-1-1-2}. 
Here $\Omega$ is defined as a square $(1.5,2.5) \times (1.0,2.0)$ with three removed and four filled circular holes. It is worth noting that there are two circles that are nearly touching at the red point $(1.935,1.78)$. This poses a key difficulty, as observed in COMSOL software \cite{comsol} with a failure for mesh generation.
To account for realistic scenarios, we consider a case where an exact solution is not available. We set $f(\boldsymbol{x})=1$, $g(\boldsymbol{x})=0$ and $a(\boldsymbol{x})= 1.0\mathbb{I}_{\Omega_{1}} + 2.0\mathbb{I}_{\Omega_{2}}$ in equation \eqref{pde-elliptic}. The obtained numerical solution and its first-order derivatives using RFM are illustrated in Figure \ref{figure3-1-1-3}.

As recorded in Table \ref{table3-1-1-2}, RFM shows a clear trend of numerical convergence. We use the RFM solution in the case with largest parameters as the reference solution and the error is about $1\%$.
\begin{table}[htbp]
	\caption{\label{table3-1-1-2} RFM results for the two-dimensional elliptic equation with complex interface geometry.} \centering
	\begin{small}
		\begin{tabular}{|c|c|ccc|}
			\hline
			$M$ & $N$ & $u$ error & $u_{x}$ error & $u_{y}$ error\\
			\hline
			\multirow{3}{*}{51200} 
			& 184596  &  4.39E-1  &  7.62E-1  &  7.64E-1 \\
			& 259920  &  2.44E-2  &  4.18E-2  &  4.05E-2 \\
			& 348038  &  7.97E-3  &  2.04E-2  &  1.62E-2 \\
			& 448966  &  \multicolumn{3}{c|}{Reference solution} \\
			\hline
		\end{tabular}
	\end{small}
\end{table}

\subsubsection{Stokes interface problem with $h_1(\boldsymbol{x})=0$ and $h_2(\boldsymbol{x})\neq0$}
\label{sec3-1-2}
The Stokes interface problem arises from multi-phase incompressible flow with density and viscosity variations across the interface between two different fluids. Consider two-phase Stokes flow equation defined by
\begin{equation}
	\begin{cases}
		\begin{aligned}
			-\mu(\boldsymbol{x}) \Delta \boldsymbol{u}(\boldsymbol{x})+\nabla p(\boldsymbol{x}) & =\boldsymbol{f}_{i}(\boldsymbol{x}), & & \boldsymbol{x}\in\Omega_{i}, \\
			\nabla \cdot \boldsymbol{u}(\boldsymbol{x}) & =0, & & \boldsymbol{x}\in \Omega_{1} \cup \Omega_{2}, \\
			\llbracket \boldsymbol{u}(\boldsymbol{x}) \rrbracket & =\boldsymbol{0}, & & \boldsymbol{x}\in \Gamma, \\
			\llbracket \boldsymbol{\sigma}(\boldsymbol{u}(\boldsymbol{x}), p(\boldsymbol{x})) \boldsymbol{n}(\boldsymbol{x}) \rrbracket & = \boldsymbol{h}_2(\boldsymbol{x}), & & \boldsymbol{x}\in \Gamma, \\
			\boldsymbol{u}(\boldsymbol{x}) & = \boldsymbol{g}(\boldsymbol{x}), & & \boldsymbol{x}\in \partial \Omega .
		\end{aligned}
	\end{cases}
	\label{pde-stokes}
\end{equation}
Here, $\boldsymbol{u}, p$ and $\boldsymbol{f}_i$ represent velocity, pressure, and external force, respectively. The stress tensor is defined as $\boldsymbol{\sigma}(\boldsymbol{u}, p)=-p \boldsymbol{I}+\mu(\nabla \boldsymbol{u}+(\nabla \boldsymbol{u})^T)$, and the viscosity $\mu$ is assumed to be a piecewise constant across the interface
\begin{equation*}
	\mu(\boldsymbol{x})= \begin{cases}\mu_{1}, & \boldsymbol{x} \in \Omega_{1} \\ \mu_{2}, & \boldsymbol{x} \in \Omega_{2}\end{cases}.
\end{equation*}
Unlike the examples in Section \ref{sec3-1-1}, Stokes interface problem \eqref{pde-stokes} contains a first-order jump condition across the interface, i.e., $h_2(\boldsymbol{x})\neq0$.

For the first example, we consider a two-dimensional Stokes interface problem in literature \cite{DBLP:journals/corr/abs-2302-08022} and use three different cases to demonstrate the robustness of our algorithm with increasing magnitude of heterogeneity:
\begin{itemize}
	\item[I] $\quad \mu_{1}=1, \quad \mu_{2}=10$,
	
	\item[II] $\quad \mu_{1}=1, \quad \mu_{2}=100$,
	
	\item[III] $\quad \mu_{1}=1, \quad \mu_{2}=1000$.
\end{itemize}

We consider a circular interface $F(\boldsymbol{x})=x^2+y^2-r^2=0$ with a radius of $r=1$, located at the center of the square domain $\Omega=(-2.0,2.0)^2$. The exact velocity and pressure are given by
\begin{equation*}
	\begin{aligned}
		u(\boldsymbol{x})=&\begin{cases}\frac{y}{4}(x^2+y^2), & \boldsymbol{x}\in\Omega_{1}\\
			\frac{y}{r}-\frac{3 y}{4}, & \boldsymbol{x}\in\Omega_{2}\end{cases},\\
		v(\boldsymbol{x})=&\begin{cases}-\frac{x y^2}{4}, & \boldsymbol{x}\in\Omega_{1}\\
			-\frac{x}{r}+\frac{x}{4}(3+x^2), & \boldsymbol{x}\in\Omega_{2}\end{cases},\\
		p(\boldsymbol{x})=&\begin{cases}5.0, & \boldsymbol{x}\in\Omega_{1} \\
			 (-\frac{3}{4} x^3+\frac{3}{8} x) y, & \boldsymbol{x}\in\Omega_{2}\end{cases}.
	\end{aligned}
\end{equation*}

\begin{figure}[htbp]
	\centering
	\subfigure[$u$]{
		\includegraphics[width=0.28\textwidth]{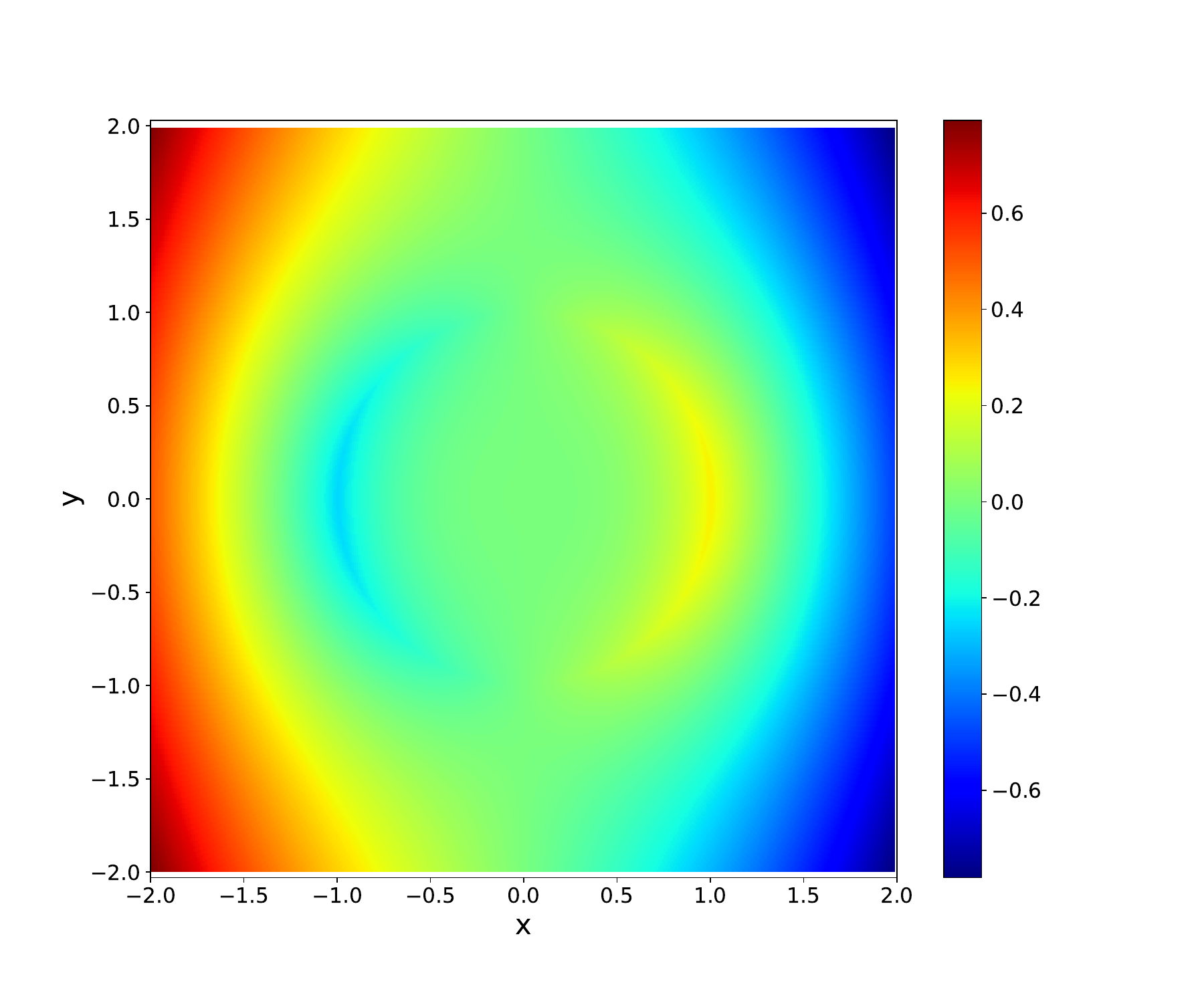}
	}
	\quad
	\subfigure[$v$]{
		\includegraphics[width=0.28\textwidth]{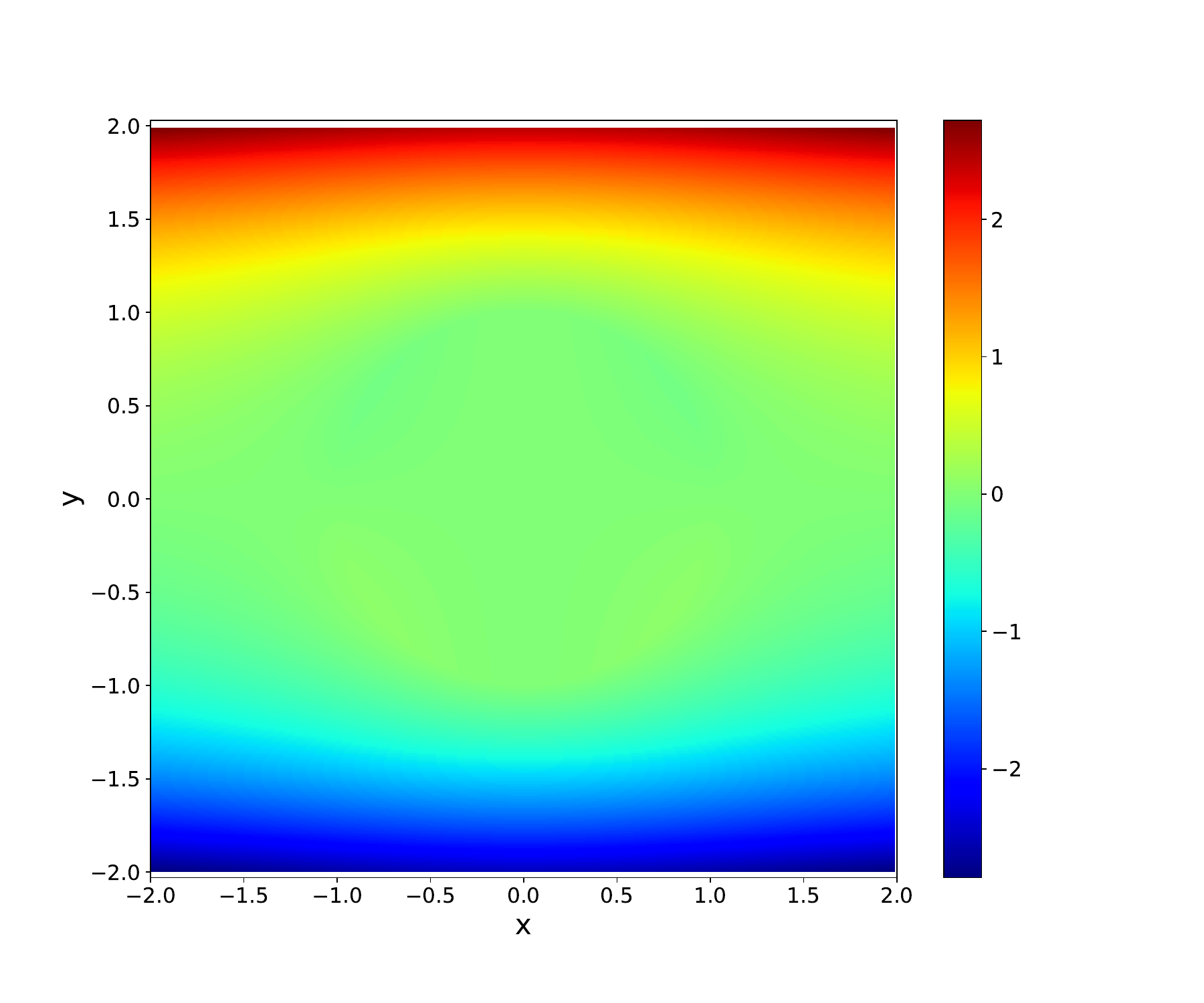}
	}
	\quad
	\subfigure[$p$]{
		\includegraphics[width=0.28\textwidth]{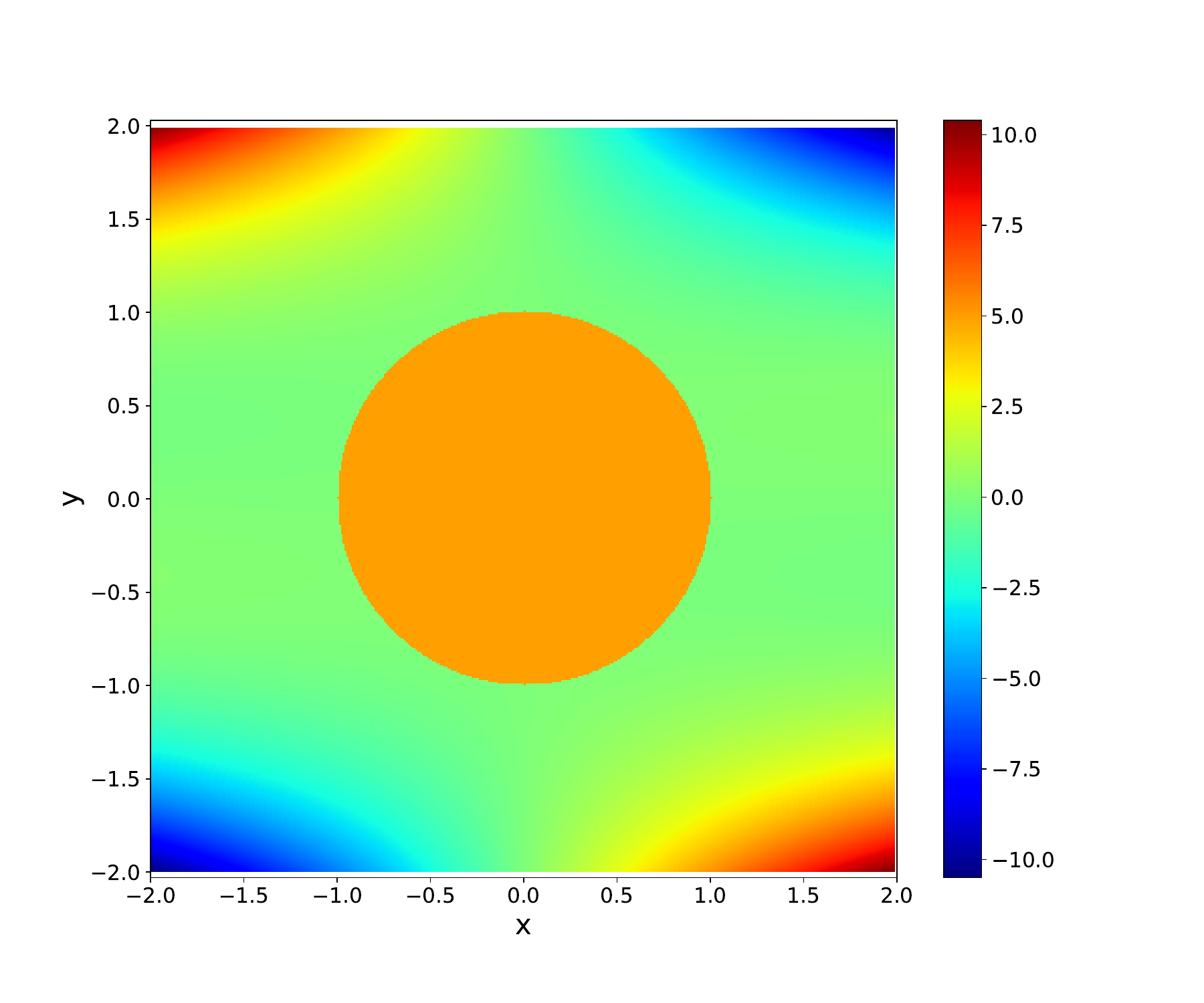}
	}
	\caption{Numerical solution for the Stokes interface problem.}
	\label{figure3-1-2-1}
\end{figure}

\begin{table}[htbp]
	\caption{\label{table3-1-2-1} Comparison of RFM and KFBI for the Stokes interface problem.} \centering
	\begin{small}
		\begin{tabular}{|c|c|cc|ccc|}
			\hline
			Case & Method & $M$ & $N$ & $u$ error  & $v$ error & $p$ error \\
			
			\hline
			\multirow{9}{*}{I} &\multirow{4}{*}{RFM} 
			& 9600 & 6801  &  1.12E-2  & 4.32E-3  & 1.05E-1   \\
			& & 38400 & 6801  &  5.96E-5  & 1.61E-5  & 6.59E-4  \\
			& & 9600 & 84801 &  6.92E-6  & 1.71E-6  & 3.33E-4   \\
			& & 38400 & 84801 &  1.18E-8  & 3.30E-9  & 4.84E-8   \\
			
			\cline{2-7}
			&\multirow{5}{*}{KFBI} 
			& 49152    & 49152    &  \multicolumn{2}{c}{3.30E-4}  & 2.63E-2   \\
			& & 196608   & 196608   &  \multicolumn{2}{c}{8.28E-5}  & 1.31E-2  \\
			& & 786432   & 786432   &  \multicolumn{2}{c}{1.96E-5}  & 6.51E-3   \\
			& & 3145728  & 3145728  &  \multicolumn{2}{c}{4.75E-6}  & 3.25E-3   \\
			& & 12582912 & 12582912 &  \multicolumn{2}{c}{1.15E-6}  & 1.62E-3   \\

			\hline
			\multirow{9}{*}{II} &\multirow{4}{*}{RFM} 
			& 9600 & 6801  &  9.01E-3  & 1.96E-3  & 4.15E-1   \\
			& & 38400 & 6801  &  1.25E-4  & 3.85E-5  & 8.98E-3  \\
			& & 9600 & 84801 &  4.93E-6  & 1.48E-6  & 2.63E-3   \\
			& & 38400 & 84801 &  1.52E-8  & 3.21E-9  & 5.99E-7   \\
			
			\cline{2-7}
			&\multirow{5}{*}{KFBI} 
			& 49152    & 49152    &  \multicolumn{2}{c}{3.47E-4}  & 2.63E-1   \\
			& & 196608   & 196608   &  \multicolumn{2}{c}{8.69E-5}  & 1.31E-1  \\
			& & 786432   & 786432   &  \multicolumn{2}{c}{2.06E-5}  & 6.51E-2   \\
			& & 3145728  & 3145728  &  \multicolumn{2}{c}{4.98E-6}  & 3.25E-2   \\
			& & 12582912 & 12582912 &  \multicolumn{2}{c}{1.20E-6}  & 1.62E-2   \\

			\hline
			\multirow{9}{*}{III} &\multirow{4}{*}{RFM} 
			& 9600 & 6801  &  7.03E-3  & 1.85E-3  & 2.33E+0   \\
			& & 38400 & 6801  &  7.33E-5  & 2.25E-5  & 4.12E-2  \\
			& & 9600 & 84801 &  4.49E-6  & 1.44E-6  & 3.31E-2   \\
			& & 38400 & 84801 &  5.97E-9  & 1.64E-9  & 2.76E-6   \\
			
			\cline{2-7}
			&\multirow{5}{*}{KFBI} 
			& 49152    & 49152    &  \multicolumn{2}{c}{3.49E-4}  & 2.63E+0   \\
			& & 196608   & 196608   &  \multicolumn{2}{c}{8.73E-5}  & 1.31E+0  \\
			& & 786432   & 786432   &  \multicolumn{2}{c}{2.07E-5}  & 6.51E-1   \\
			& & 3145728  & 3145728  &  \multicolumn{2}{c}{5.00E-6}  & 3.25E-1   \\
			& & 12582912 & 12582912 &  \multicolumn{2}{c}{1.21E-6}  & 1.62E-1   \\
			\hline
		\end{tabular}
	\end{small}
\end{table}
Convergence of RFM is recorded in Table \ref{table3-1-2-1} and results of KFBI are also listed for comparison. In both methods, the velocity remains highly accurate but the pressure accuracy deteriorates significantly as the jump of viscosity coefficients increases. The accuracy of RFM outperforms that of KFBI by $2\sim3$ orders of magnitude for velocity and by $4\sim5$ orders of magnitude for pressure, respectively. Moreover, for the same accuracy requirement, the degrees of freedom in RFM is $3\sim4$ orders of magnitude smaller than KFBI method.

Next, we investigate whether the smoothness of the solution near the interface affects the accuracy of RFM in a three-dimensional Stokes interface problem defined on the cubic domain $\Omega = (-1.0,1.0)^3$ with a spherical interface $F(\boldsymbol{x})=x^2+y^2+z^2-1=0$. The viscosity is set as $\mu_1=1,\mu_2=10$.

Two exact solutions are considered for comparison: a smooth solution $\boldsymbol{u}_{s}=(u_{0}, v_{0}, w_{0})^T$ and a non-smooth one $\boldsymbol{u}_{n}=(u, v, w)^T$
\begin{equation*}
	\begin{aligned}
		u_{0}(\boldsymbol{x})=&x^2\left(2-x^2\right)+4 x y\left(x^2+y^2-1\right)+z^2\left(3 z^2-6 x^2-2\right), \\
		v_{0}(\boldsymbol{x})=&x^2\left(3 x^2-6 y^2-1\right)+y^2\left(1-y^2\right), \\
		w_{0}(\boldsymbol{x})=&4 z x\left(z^2+x^2-1\right), \\ 
		u(\boldsymbol{x})=&\begin{cases}u_{0}+(y-z)(x^2+y^2+z^2-1), & \boldsymbol{x}\in\Omega_{1} \\
			u_{0}, & \boldsymbol{x}\in\Omega_{2} \end{cases}, \\
		v(\boldsymbol{x})=&\begin{cases}v_{0}+(z-x)(x^2+y^2+z^2-1), & \boldsymbol{x}\in\Omega_{1}\\
			v_{0}, & \boldsymbol{x}\in\Omega_{2} \end{cases},\\
		w(\boldsymbol{x})=&\begin{cases}w_{0}+(x-y)(x^2+y^2+z^2-1), & \boldsymbol{x}\in\Omega_{1} \\
			w_{0}, & \boldsymbol{x}\in\Omega_{2}\end{cases},\\				
		p(\boldsymbol{x})=&8 y\left(3 x^2-y^2\right)+8 x\left(3 z^2-x^2\right). \\			
	\end{aligned}
\end{equation*}	

Table \ref{table3-1-2-2} presents the influence of the solution smoothness on the convergence of RFM. It is clear that the smoothness of the exact solution has a minimal effect on the performance of RFM.
\begin{table}[htbp]
	\caption{\label{table3-1-2-2} Results of RFM for the three-dimensional Stokes interface problem with smooth and nonsmooth solutions.} \centering
	\begin{small}
		\begin{tabular}{|c|ccc|}
			\hline
			Exact solution & $u$ error  & $v$ error & $w$ error \\
			\hline
			$\boldsymbol{u}_s$ &  1.14E-5  & 1.11E-5  & 1.97E-5   \\
			$\boldsymbol{u}_n$ &  1.11E-5  & 1.55E-5  & 1.85E-5  \\
			\hline
		\end{tabular}
	\end{small}
\end{table}

For the Stokes problems on three-dimensional irregular domains, the KFBI method solves equivalent but simpler interface problems within an extended cubic region. In contrast, RFM can directly solve these problems without transforming them into equivalent interface problems. Specifically, RFM utilizes the level-set method for boundary point sampling and follows the algorithm described in \cite{JML-1-268}.
	

The Stokes problem is defined over an irregular torus domain $\Omega=F(\boldsymbol{x})<0$, where $F(\boldsymbol{x})=(c-\sqrt{(x^2+y^2)})^2+z^2-a^2$ with $a = 0.35$ and $c = 0.7$.
The exact solution is specified as
\begin{equation*}
	\begin{aligned}
		u(\boldsymbol{x})  =&-4 x y(1-x^2-y^2)-x^2(x^2+6 z^2-2)+z^2(3 z^2-2) \\
		& +\exp (\cos y)+\exp (\sin z), \\
		v(\boldsymbol{x})  =&x^2(3 x^2-6 y^2-2)-y^2(y^2-2)+\exp (\sin x), \\
		w(\boldsymbol{x})  =&-4(1-x^2-z^2) x z+\exp (\cos x), \\
		p(\boldsymbol{x})  =&\exp (1-y^2-z^3) \sin (x^2+1) . \\			
	\end{aligned}
\end{equation*}


\begin{table}[htbp]
	\caption{\label{table3-1-2-3} Comparison of RFM and KFBI for three-dimensional Stokes problem.} \centering
	\begin{small}
		\begin{tabular}{|c|cc|ccc|}
			\hline
			Method & $M$ & $N$ & $u$ error  & $v$ error & $w$ error \\
			
			\hline
			\multirow{3}{*}{RFM} 
			& 25600 & 129088  &  5.34E-6  & 1.53E-5  & 2.97E-6   \\
			& 51200 & 129088  &  2.57E-6  & 7.26E-6  & 1.22E-6  \\
			& 102400& 129088  &  1.47E-6  & 4.27E-6  & 1.04E-6   \\
			
			\hline
			\multirow{3}{*}{KFBI} 
			& 8388608    & 8388608    & 2.22E-5 & 7.64E-5 & 2.62E-5   \\
			& 67108864   & 67108864   & 4.74E-6 & 1.26E-5 & 3.95E-6  \\
			& 536870912  & 536870912  & 1.03E-6 & 2.50E-6 & 9.14E-7   \\
			\hline
		\end{tabular}
	\end{small}
\end{table}
Table \ref{table3-1-2-3} records the error in RFM and KFBI. It is evident that RFM achieves the same accuracy with significantly fewer degrees of freedom, by three orders of magnitude.

\subsubsection{Elasticity interface problem with $h_1(\boldsymbol{x})\neq0$ and $h_2(\boldsymbol{x})\neq0$}
\label{sec3-1-3}
Elastic interface problems have a wide range of science and engineering applications, for example, the dynamics of crystalline materials, the simulation of microstructural evolution, and the modelling of atomic interactions. In this section, we consider the following elasticity interface problem \eqref{pde-elasticity}
\begin{equation}
	\begin{cases}
		\begin{aligned}
			-\nabla \cdot \sigma(\boldsymbol{u}(\boldsymbol{x})) & =\boldsymbol{f}_i(\boldsymbol{x}), & & \boldsymbol{x}\in \Omega_{i} , \\
			\llbracket u(\boldsymbol{x}) \rrbracket&=\boldsymbol{h}_1(\boldsymbol{x}),& & \boldsymbol{x}\in \Gamma,\\
			\llbracket \sigma(\boldsymbol{u}) \boldsymbol{n} \rrbracket&=\boldsymbol{h}_2(\boldsymbol{x}), & & \boldsymbol{x}\in \Gamma, \\
			u(\boldsymbol{x})&=g(\boldsymbol{x}), & & \boldsymbol{x}\in  \partial \Omega,
		\end{aligned}
	\end{cases}
	\label{pde-elasticity}
\end{equation}
where $\boldsymbol{u}$ and $\boldsymbol{f}_i$ represent the displacement and body force, respectively. The strain tensor and the stress tensor are defined as $\epsilon(\boldsymbol{u})=\frac{1}{2}(\nabla \boldsymbol{u}+\nabla \boldsymbol{u}^{\top})$ and $\sigma(\boldsymbol{u})=\lambda(\nabla \cdot \boldsymbol{u}) \mathbf{I}+2 \mu \epsilon(\boldsymbol{u})$, where $\lambda$ and $\mu$ are the Lam\'e parameters assumed to be piecewise constants
\begin{equation*}
\lambda=\begin{cases}
	\lambda_1, & \boldsymbol{x}\in \Omega_1\\
	\lambda_2, & \boldsymbol{x}\in \Omega_2\end{cases},
\quad \mu= \begin{cases}\mu_1, & \boldsymbol{x}\in \Omega_1\\
	\mu_2, & \boldsymbol{x}\in \Omega_2\end{cases}.
\end{equation*}

We adopt the same condition as in \cite{2022High}, where both interface conditions $\boldsymbol{h}_1(\boldsymbol{x})$ and $\boldsymbol{h}_2(\boldsymbol{x})$ are non-zero. Let the interface $F(\boldsymbol{x})=0$ be a ball centered at $(0.5,0.5,0.5)^T$ with radius $r=0.25$, and the entire domain $\Omega=(0.0,1.0)^3$. The exact solution is
\begin{equation*}
	\begin{aligned}
		u(\boldsymbol{x})=&\begin{cases}-\cos (x^2) \mathrm{e}^{-y^2} \sin (2 \pi z), & \boldsymbol{x}\in\Omega_{1} \\
			(-\sin (x^2) \mathrm{e}^{y^2} \cos (2 \pi z), & \boldsymbol{x}\in\Omega_{2}\end{cases}, \\
		v(\boldsymbol{x})=&\begin{cases}-\cos (y^2) \mathrm{e}^{-x^2} \sin (2 \pi z), & \boldsymbol{x}\in\Omega_{1} \\
			-\sin (y^2) \mathrm{e}^{x^2} \cos (2 \pi z), & \boldsymbol{x}\in\Omega_{2}\end{cases},\\
		w(\boldsymbol{x})=&\begin{cases}\cos (y^2) \mathrm{e}^{-z^2} \sin (2 \pi x), & \boldsymbol{x}\in\Omega_{1} \\
			\sin (y^2) \mathrm{e}^{z^2} \cos (2 \pi x), & \boldsymbol{x}\in\Omega_{2}\end{cases}.\\
	\end{aligned}
\end{equation*}
The Lam\'e parameters are chosen as $\left(\lambda_1, \lambda_2\right)=\left(\mu_1, \mu_2\right)=(1,100)$.

Results of RFM and penalty FEM \cite{2022High} are shown in Table \ref{table3-1-3}.
\begin{table}[htbp]
	\caption{\label{table3-1-3} Comparison of RFM and FEM for the three-dimensional elasticity interface problem.} \centering
	\begin{small}
		\begin{tabular}{|c|cc|ccc|}
			\hline
			Method & $M$ & $N$ & $u$ error  & $v$ error & $w$ error \\
			
			\hline
			\multirow{2}{*}{RFM} 
			& 38400 & 344820  &  1.91E-6  & 1.93E-6  & 1.79E-6  \\
			& 57600 & 344820  &  1.52E-7  & 1.19E-7  & 1.72E-7  \\
			
			\hline
			\multirow{2}{*}{penalty FEM} 
			& 3220614  & 3220614  & \multicolumn{3}{c|}{2.34E-5}  \\
			& 1369446  & 1369446  & \multicolumn{3}{c|}{4.76E-6}  \\
			\hline
		\end{tabular}
	\end{small}
\end{table}
It is demonstrated that RFM achieves high accuracy even for interface problems with discontinuous solutions, i.e., $\boldsymbol{h}_1(\boldsymbol{x})\neq 0$, and the reduction in terms of degrees of freedom is also evident in this example.

\subsection{Time-dependent Interface Problems}
\label{sec3-2}

Time-dependent interface problems often involve interfaces that evolve over time. These types of problems have numerous applications in medicine and engineering, such as blood flow dynamics and free surface phenomena. Moving interface problems poses challenges in mesh generation. It is expected that RFM has a clear advantage in this case.

To illustrate the evolution behavior, we consider a fixed domain $\Omega \subset \mathbb{R}^2$ divided into two time-varying subdomains $\Omega_1(t)$ and $\Omega_2(t)$, which are separated by an evolving interface $\Gamma(t)$.
Suppose there is a certain advection velocity $\boldsymbol{w}(\boldsymbol{X},t)$ that drives the evolution of $\Gamma(t)$, namely,
\begin{equation*}
	\frac{\mathrm{d}\boldsymbol{X}}{\mathrm{d} t}=\boldsymbol{w}(\boldsymbol{X}, t), \quad \boldsymbol{X}\in \Gamma(t).
	\label{evaluation}
\end{equation*}

This section includes three experiments: moving interface problem with decoupled advection velocities $\boldsymbol{w}$ and interface topological change; dynamic interface problem with large deformation; linear fluid-solid interaction problem with a complex geometry. For the second example, we compare RFM with the interface tracking algorithm. Hyper-parameters of RFM are given in Appendix \ref{sec::appendix} for time-dependent problems.

\subsubsection{Moving interface problem with topological change}
\label{sec3-2-1}
Consider the following parabolic moving interface model where the advection velocity $\boldsymbol{w}$ is independent of the equation
\begin{equation}
	\begin{cases}
		\begin{aligned}
			\partial_t u(\boldsymbol{x},t)-\nabla \cdot(\beta_{i} \nabla u(\boldsymbol{x},t))&=f_{i}(\boldsymbol{x},t), &&\quad  \boldsymbol{x} \in  \Omega_{i}(t), \quad t \in [0, T],\\
			\llbracket u(\boldsymbol{x},t) \rrbracket&=h_1(\boldsymbol{x},t), &&\quad \boldsymbol{x} \in \Gamma(t), \quad t \in [0, T],\\
			\llbracket \beta \nabla u(\boldsymbol{x},t) \cdot \boldsymbol{n}(\boldsymbol{x},t) \rrbracket&=h_2(\boldsymbol{x},t), &&\quad \boldsymbol{x} \in \Gamma(t), \quad t \in [0, T],\\
			u(\boldsymbol{x}, 0)&=u_0(\boldsymbol{x}), && \quad \boldsymbol{x} \in \Omega,\\
			u(\boldsymbol{x}, t)&=g(\boldsymbol{x}, t), && \quad \boldsymbol{x} \in \partial \Omega, \quad t \in [0, T],
		\end{aligned}
	\end{cases}
	\label{pde-parabolic}
\end{equation}
with the constants fixed as $\beta_1=1,\beta_2=10$, and the entire domain $\Omega=(-1.0,1.0)^2$.

We begin with a simple interface evolution where the circular interface $F(\boldsymbol{x},t)=(x-0.3 \cos (\pi t))^2+(y-0.3 \sin (\pi t))^2-(\pi / 6)^2=0, t\in (0,1)$ moves with both rotational and translational motions. The immersed finite element method \cite{doi:10.1137/20M133508X} applies difference schemes in the temporal dimension and generates unfitted meshes at each time step to handle interface movement. In contrast, RFM increases the flexibility of the method by handling both spatial and temporal dimensions within the same framework, while maintaining high accuracy.


For this problem, we use the same exact solution as in \cite{doi:10.1137/20M133508X}
\begin{tiny}
	\begin{equation*}
		u(\boldsymbol{x}, t)=\begin{cases}\frac{((x-0.3 \cos (\pi t))^2+(y-0.3 \sin (\pi t))^2)^{5 / 2}(\pi / 6)^{-1}}{\beta_{1}}, & \boldsymbol{x}\in\Omega_{1}(t) \\
			\frac{\left((x-0.3 \cos (\pi t))^2+(y-0.3 \sin (\pi t))^2\right)^{5 / 2}(\pi / 6)^{-1}}{\beta_{2}}+(\pi / 6)^4(\frac{1}{\beta_{1}}-\frac{1}{\beta_{2}}), & \boldsymbol{x}\in\Omega_{2}(t)\end{cases}. \\
	\end{equation*}
\end{tiny}

\begin{figure}[htbp]
	\centering
	\subfigure[Circle interface and its movement]{
		\includegraphics[width=0.4\textwidth]{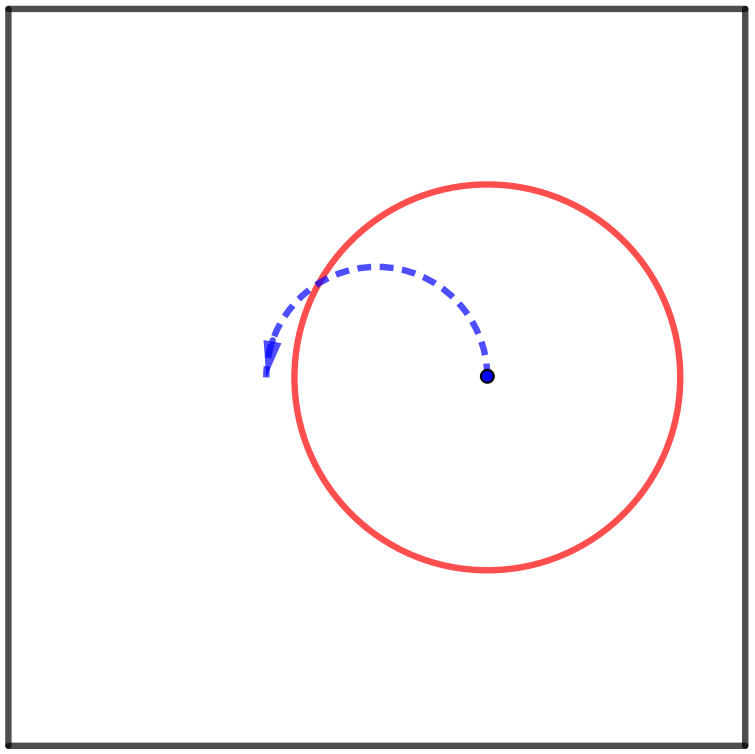}
	}
	\quad
	\subfigure[Errors in IFM]{
		\includegraphics[width=0.5\textwidth]{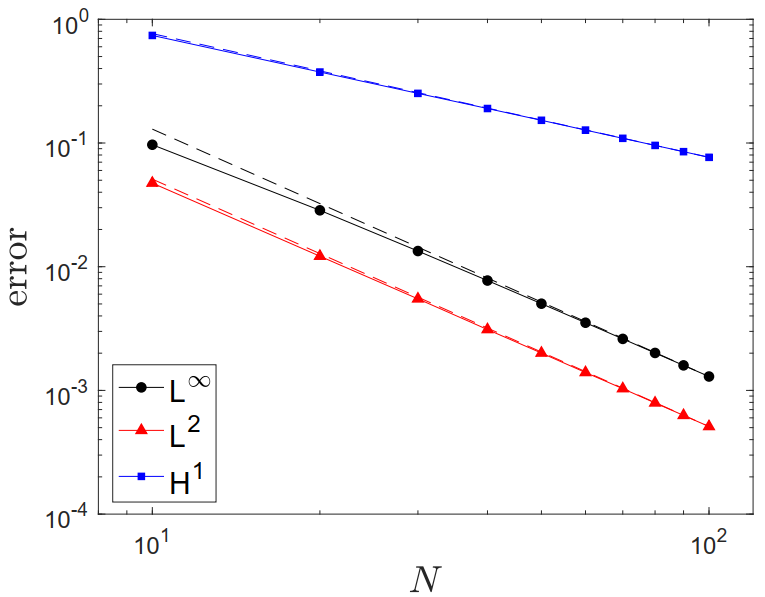}
	}
	\caption{The parabolic moving interface problem and results of the immersed finite element method.}
	\label{fig3-2-1-1}
\end{figure}

\begin{table}[htbp]
	\caption{\label{table3-2-1-1} $L^{2}$ error of RFM for the parabolic moving interface problem.} \centering
	\begin{small}
		\begin{tabular}{|c|c|ccc|}
			\hline
			$M$ & $N$ & $u$ error & $u_x$ error & $u_y$ error \\
			\hline
			12800 & 106000  &  1.96E-4   & 7.64E-4  & 8.68E-4 \\
			25600 & 375360  &  4.52E-5   & 4.04E-4  & 4.44E-4 \\
			\hline
		\end{tabular}
	\end{small}
\end{table}
Results of RFM are shown in Table \ref{table3-2-1-1}. RFM achieves an accuracy of $10^{-4}$ with $10^4$ degrees of freedom, while IFM only achieves such an accuracy with $10^8$ degrees of freedom. Furthermore, RFM has a significant advantage in fitting first-order derivatives, with almost the same magnitude of accuracy. In contrast, IFM has a drop of more than $10^{-2}$ in the $H^1$ norm.

\begin{figure}[htbp]
	\centering
	\subfigure[$t=0.0$]{
		\includegraphics[width=0.28\textwidth]{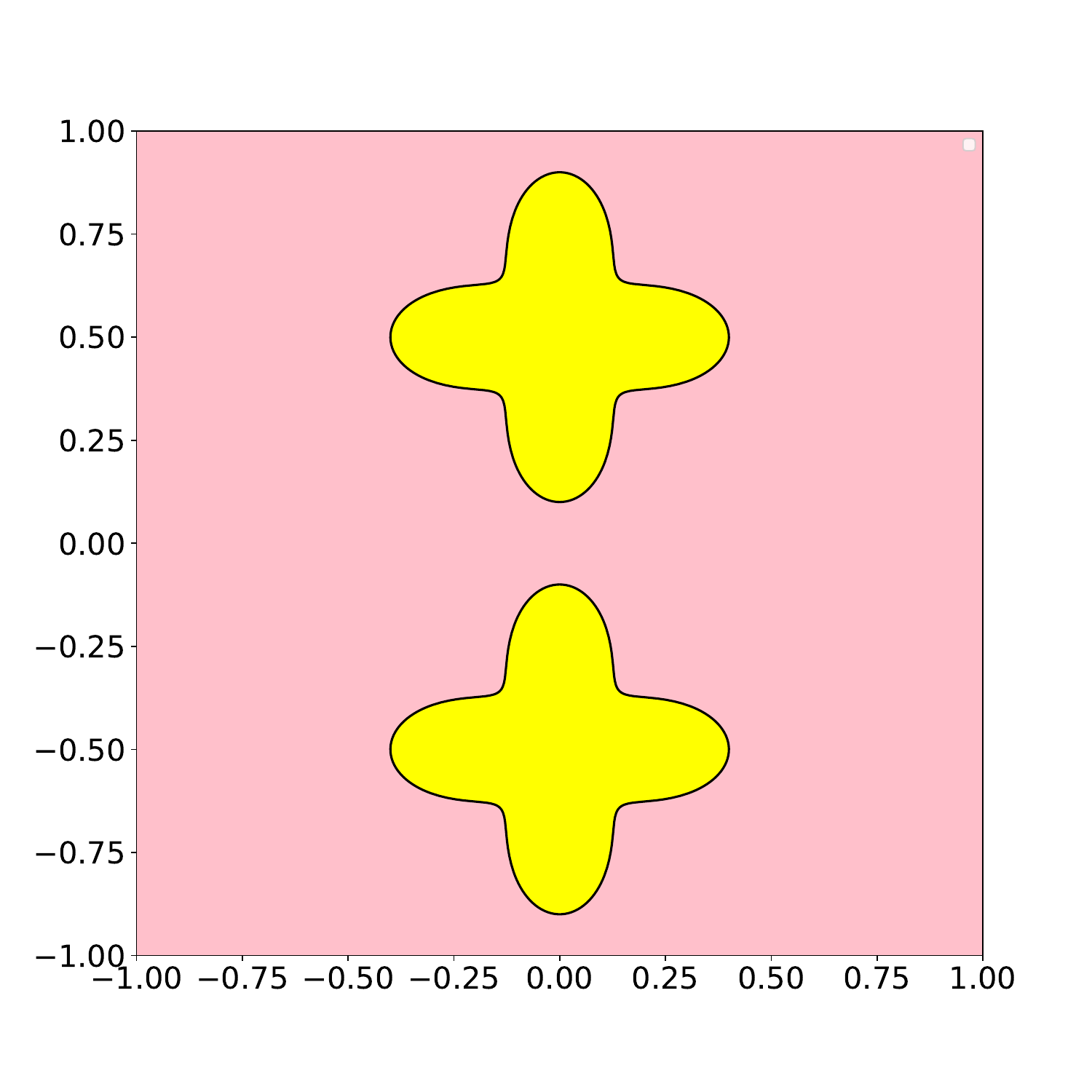}
	}
	\quad
	\subfigure[$t=0.2$]{
		\includegraphics[width=0.28\textwidth]{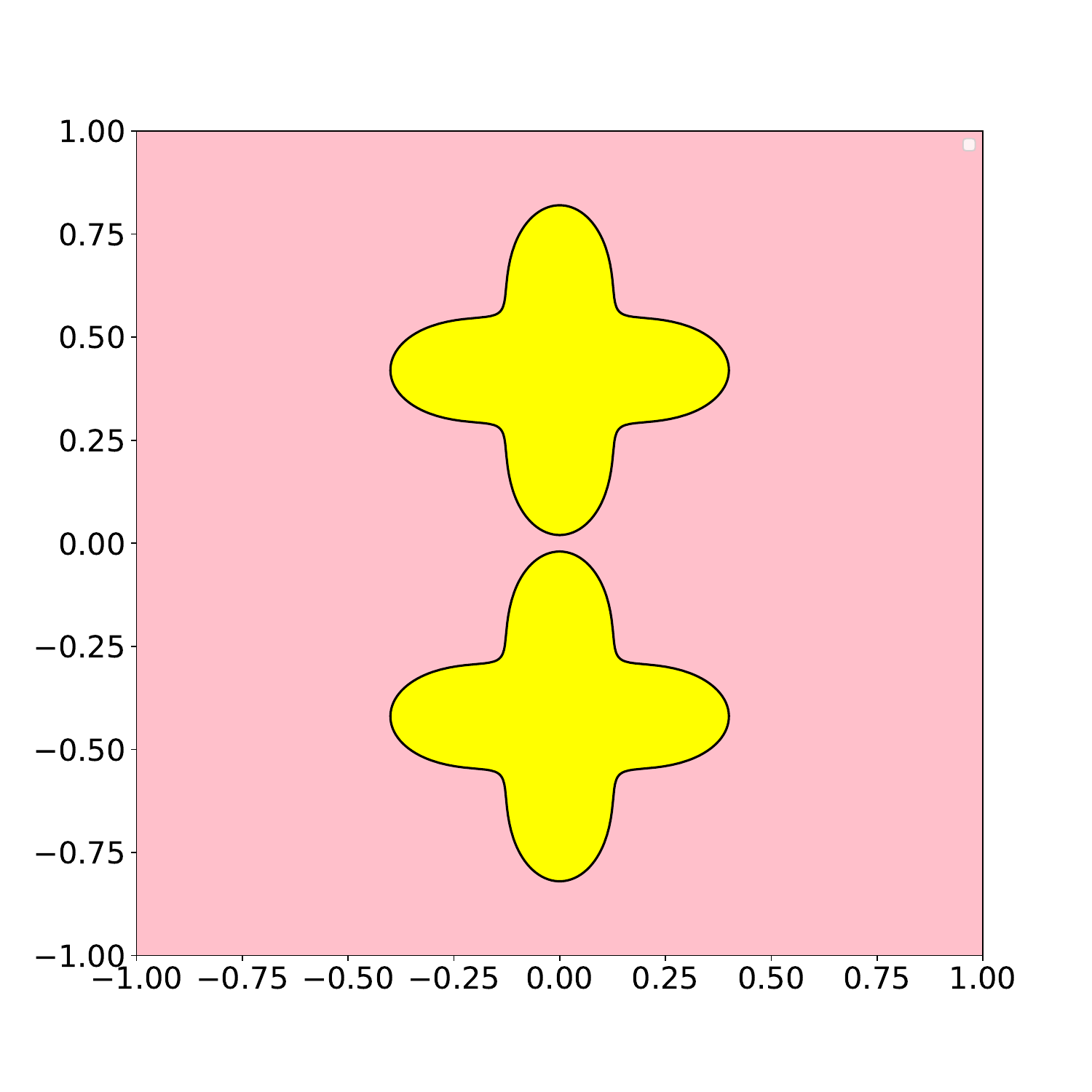}
	}
	\quad
	\subfigure[$t=0.4$]{
		\includegraphics[width=0.28\textwidth]{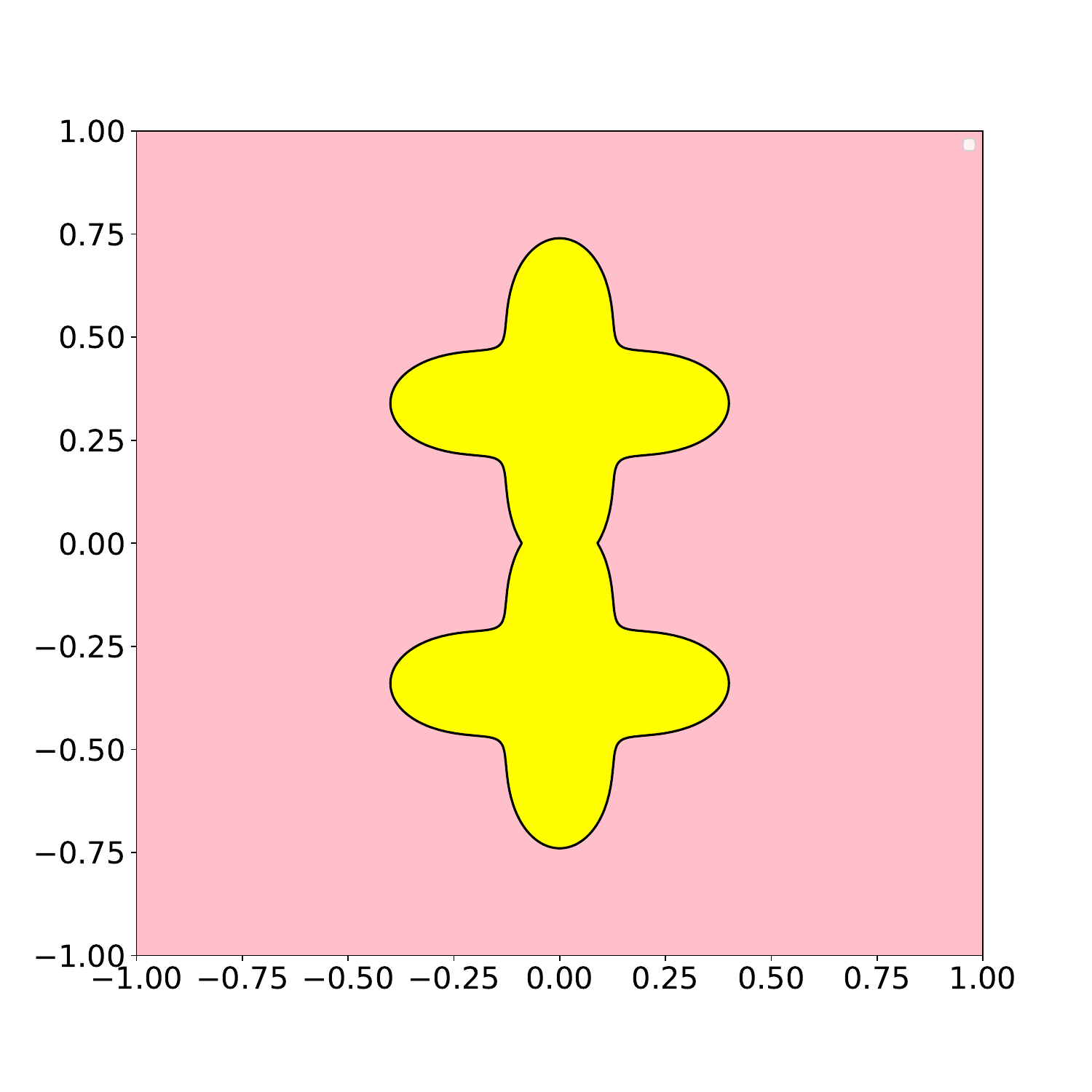}
	}
	\quad
	\subfigure[$t=0.6$]{
		\includegraphics[width=0.28\textwidth]{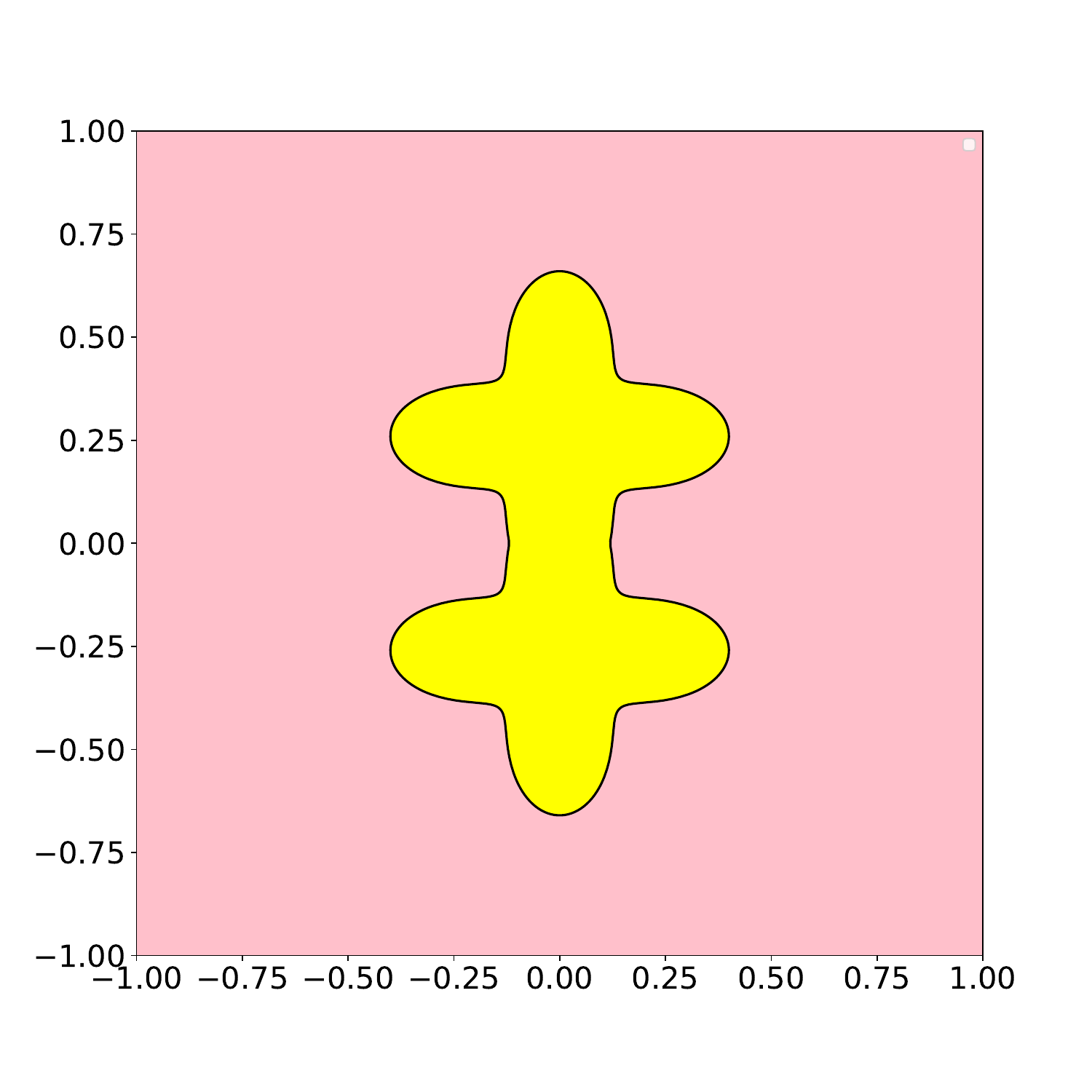}
	}
	\quad
	\subfigure[$t=0.8$]{
		\includegraphics[width=0.28\textwidth]{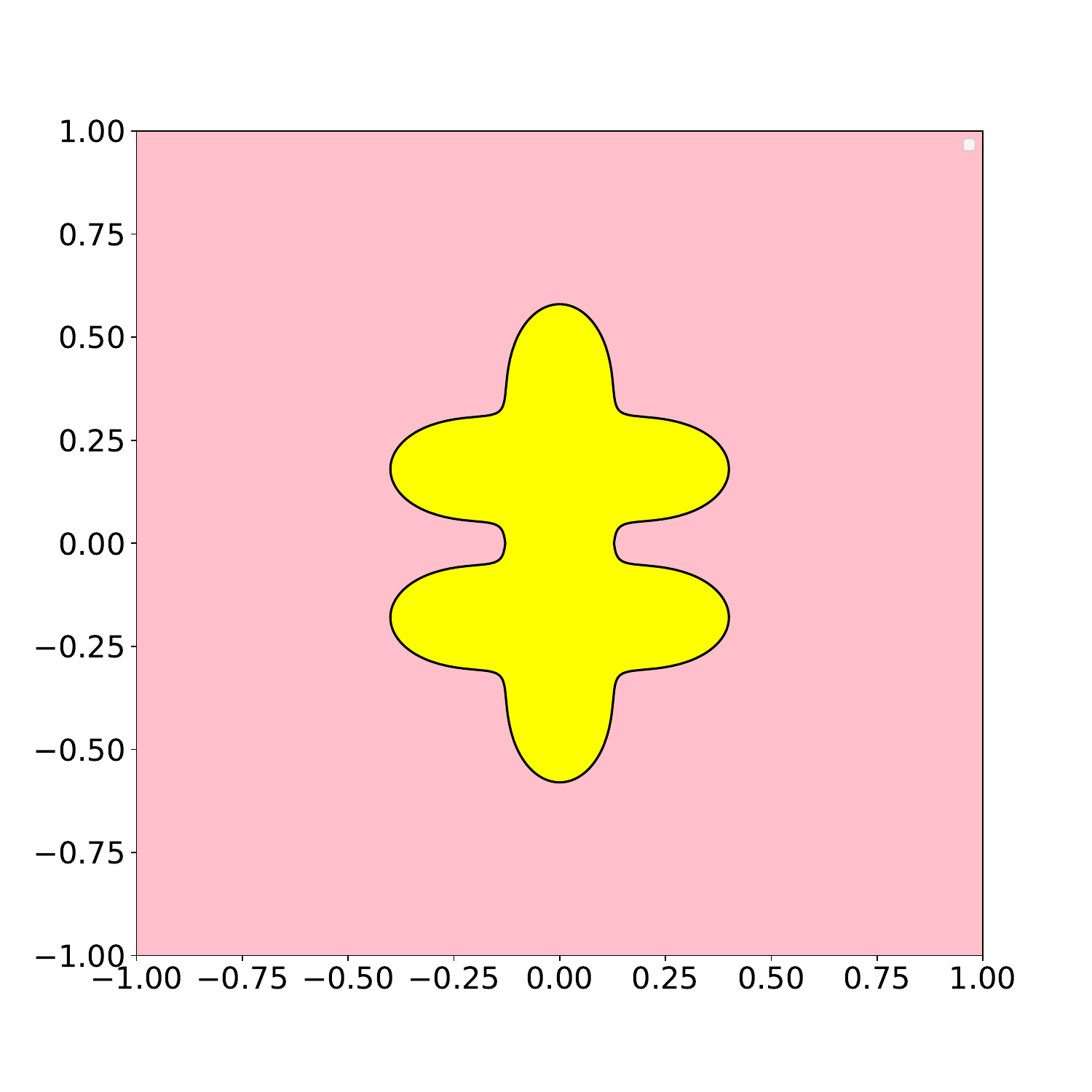}
	}
	\quad
	\subfigure[$t=1.0$]{
		\includegraphics[width=0.28\textwidth]{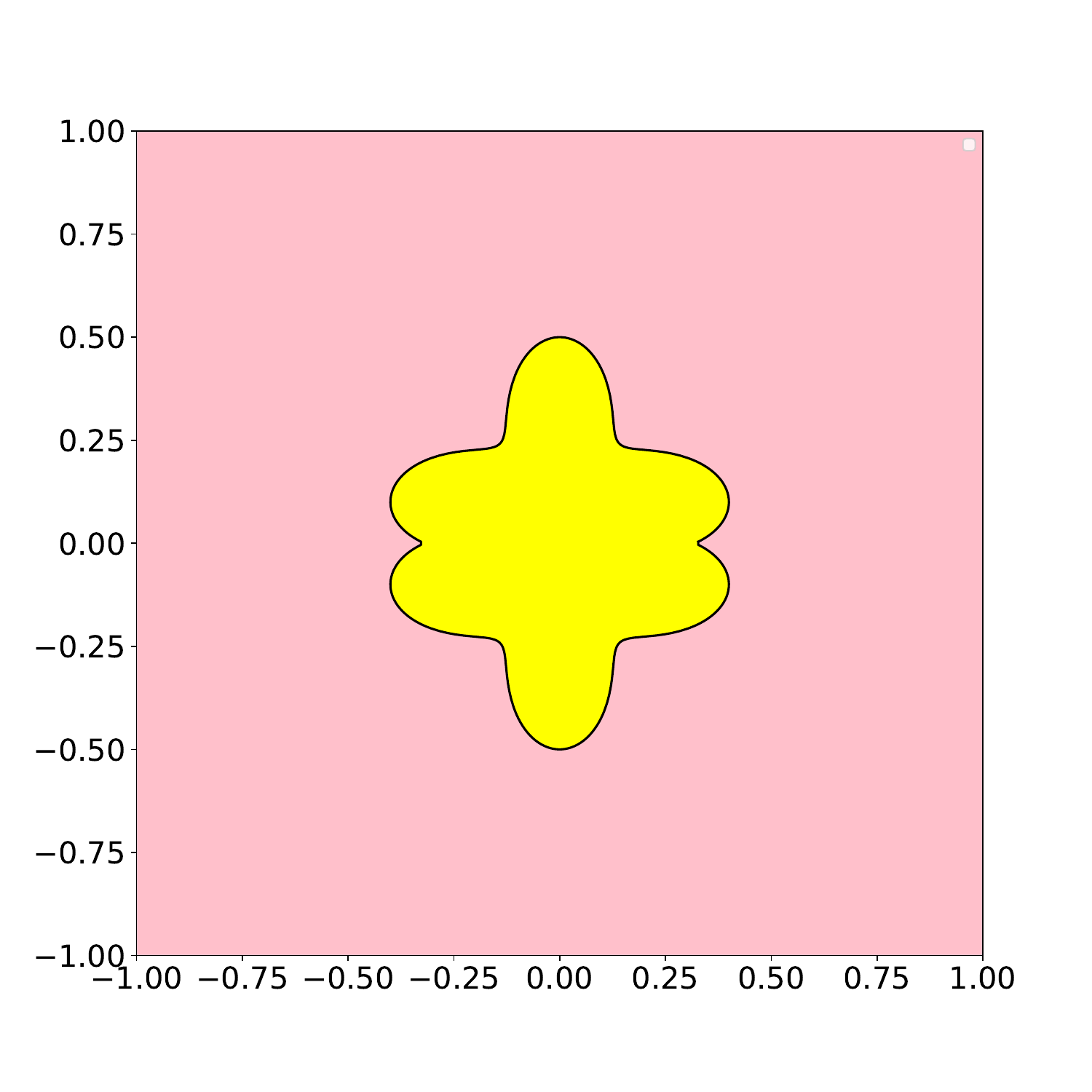}
	}
	\caption{Merging process of two interface at $t = 0.0$, $0.2$, $0.4$, $0.6$, $0.8$, and $1.0$, respectively.}
	\label{fig3-2-1-2}
\end{figure}

Next, we consider a more intricate case with the topological change of two merging quadrangular interfaces.
\begin{equation*}
	\begin{aligned}
		& F_1(\boldsymbol{x}, t) = x^2 + (y - 0.5 + 0.4t)^2 - (0.3 + 0.1\cos(4\arctan \frac{y - 0.5 + 0.4t}{x}))^2=0, \\
		& F_2(\boldsymbol{x}, t) = x^2 + (y + 0.5 - 0.4t)^2 - (0.3 + 0.1\cos(4\arctan \frac{y + 0.5 - 0.4t}{x}))^2=0.
	\end{aligned}
\end{equation*}
Figure \ref{fig3-2-1-2} visualizes the merging process at $t = 0.0$, $0.2$, $0.4$, $0.6$, $0.8$, and $1.0$, respectively.

The exact solution is denoted as
\begin{equation*}
	u(\boldsymbol{x}, t)=\begin{cases}\sin(\pi x)\sin(\pi y)\sin(t), & \boldsymbol{x}\in\Omega_{1}(t) \\
		10 - x^2 - y^2 - t^2, & \boldsymbol{x}\in\Omega_{2}(t) \end{cases}.\\
\end{equation*}

\begin{table}[htbp]
	\caption{\label{table3-2-1-2} Results of RFM for the parabolic moving interface problem with topological change.} \centering
	\begin{small}
		\begin{tabular}{|c|c|ccc|}
			\hline
			$M$ & $N$ & $u$ error & $u_x$ error & $u_y$ error \\
			\hline
			800   & 141261  &  1.98E-2   & 2.65E-1  & 2.14E-1 \\
			1600  & 141261  &  9.27E-4   & 2.67E-2  & 2.06E-2 \\
			3200  & 141261  &  4.13E-4   & 8.91E-3  & 9.16E-3 \\
			6400  & 141261  &  1.20E-5   & 3.31E-4  & 3.01E-4 \\
			12800 & 141261  &  8.87E-7   & 8.52E-6  & 1.80E-5 \\
			\hline
		\end{tabular}
	\end{small}
\end{table}
From Table \ref{table3-2-1-2}, we observe that RFM maintains the spectral accuracy even in the case of time-dependent interface problems with topological change.

\subsubsection{Dynamic interface problem with large deformation}
\label{sec3-2-2}
The linear interface problem of two-phase incompressible fluids \cite{ma2021high} is defined as Osean equation \eqref{pde-osean}
\begin{tiny}
\begin{equation}
	\begin{cases}
		\begin{aligned}
			\frac{\partial \boldsymbol{u}_i(\boldsymbol{x},t)}{\partial t}+(\boldsymbol{w}(\boldsymbol{x},t) \cdot \nabla) \boldsymbol{u}_i(\boldsymbol{x},t)-\nu_i \Delta \boldsymbol{u}_i(\boldsymbol{x},t)+\nabla p_i(\boldsymbol{x},t)&=\boldsymbol{f}_i(\boldsymbol{x},t), && \quad  \boldsymbol{x} \in  \Omega_{i}(t), \quad t \in [0, T],\\
			\operatorname{div} \boldsymbol{u}(\boldsymbol{x},t)&=0, && \quad  \boldsymbol{x} \in  \Omega ,\\
			\llbracket \boldsymbol{u}(\boldsymbol{x},t)\rrbracket&=\boldsymbol{h_1}(\boldsymbol{x},t), && \quad \boldsymbol{x} \in \Gamma(t), \quad t \in [0, T],\\
			\llbracket \nu \partial_{\boldsymbol{n}} \boldsymbol{u}(\boldsymbol{x},t)-p(\boldsymbol{x},t) \boldsymbol{n}(\boldsymbol{x},t) \rrbracket&=\boldsymbol{h_2}(\boldsymbol{x},t), && \quad \boldsymbol{x} \in \Gamma(t), \quad t \in [0, T],\\
			\boldsymbol{u}(\boldsymbol{x},0)&=\boldsymbol{u}_{0}(\boldsymbol{x}), && \quad  \boldsymbol{x} \in \Omega, \\
			\boldsymbol{u}(\boldsymbol{x},t)&=\boldsymbol{g}(\boldsymbol{x}, t), && \quad \boldsymbol{x} \in \partial \Omega, \quad t \in [0, T],
		\end{aligned}
	\end{cases}
	\label{pde-osean}
\end{equation}
\end{tiny}
where $\boldsymbol{u}_i, p_i$, and $\boldsymbol{f}_i$ represent the flow velocity, pressure, and body force, respectively.

The entire domain is $\Omega=(0.0,1.0)^2$ with a time interval of $[0, T]$, where $T=1.5$. At time $t=0$, the initial subdomain $\Omega_1(0)$ is a disk of radius 0.15 and centering at $(0.5,0.75)$. The flow velocity that drives the interface $\Gamma(t)$ is denoted as
\begin{equation*}
	\boldsymbol{w}(\boldsymbol{x}, t)=\cos (\pi t / 3)\left(\sin ^2\left(\pi x\right) \sin \left(2 \pi y\right),-\sin ^2\left(\pi y\right) \sin \left(2 \pi x\right)\right)^{T}.
\end{equation*}

\begin{figure}[htbp]
	\centering
	\subfigure[$t=0.0$]{
		\includegraphics[width=0.45\textwidth]{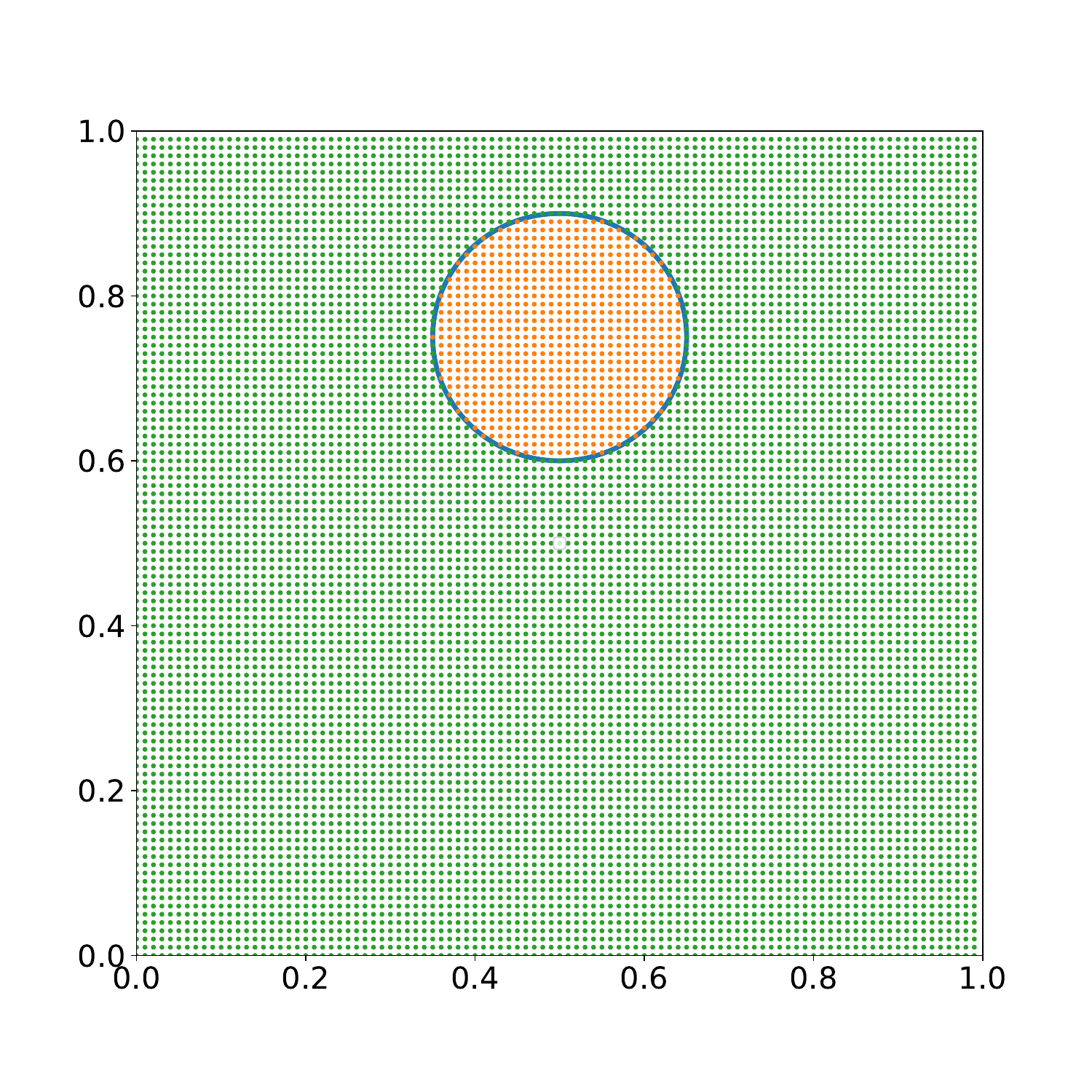}
	}
	\quad
	\subfigure[$t=0.5$]{
		\includegraphics[width=0.45\textwidth]{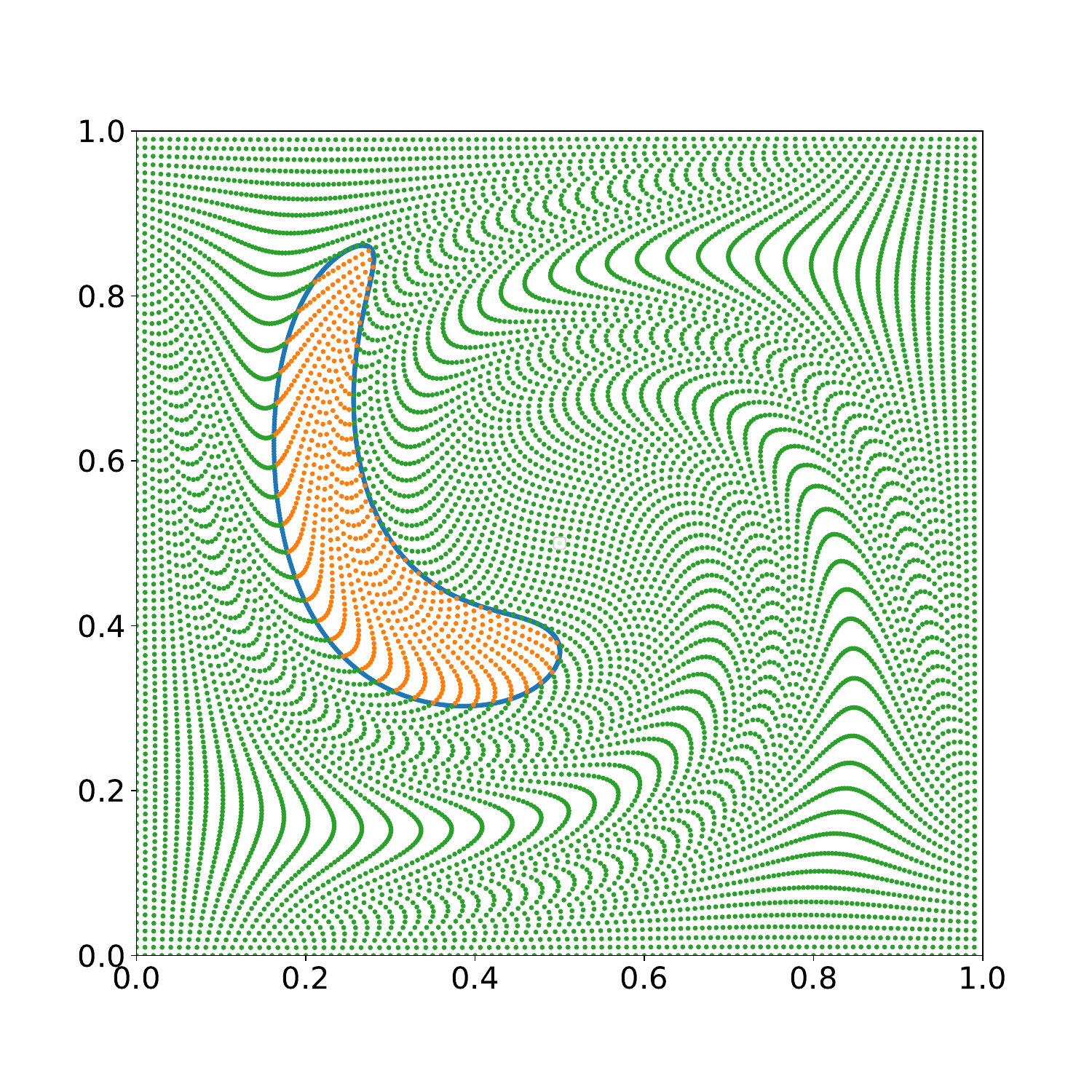}
	}
	\quad
	\subfigure[$t=1.0$]{
		\includegraphics[width=0.45\textwidth]{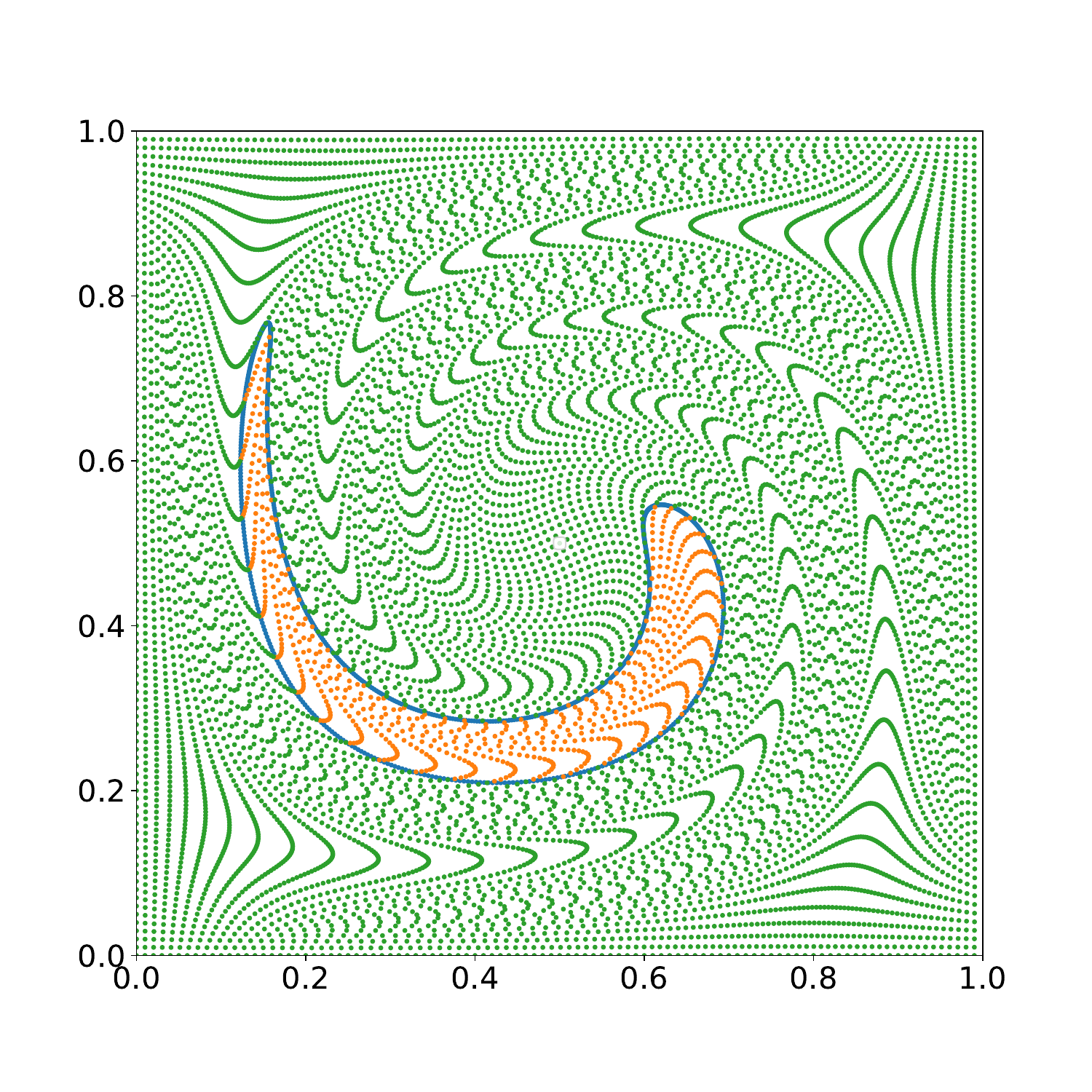}
	}
	\quad
	\subfigure[$t=1.5$]{
		\includegraphics[width=0.45\textwidth]{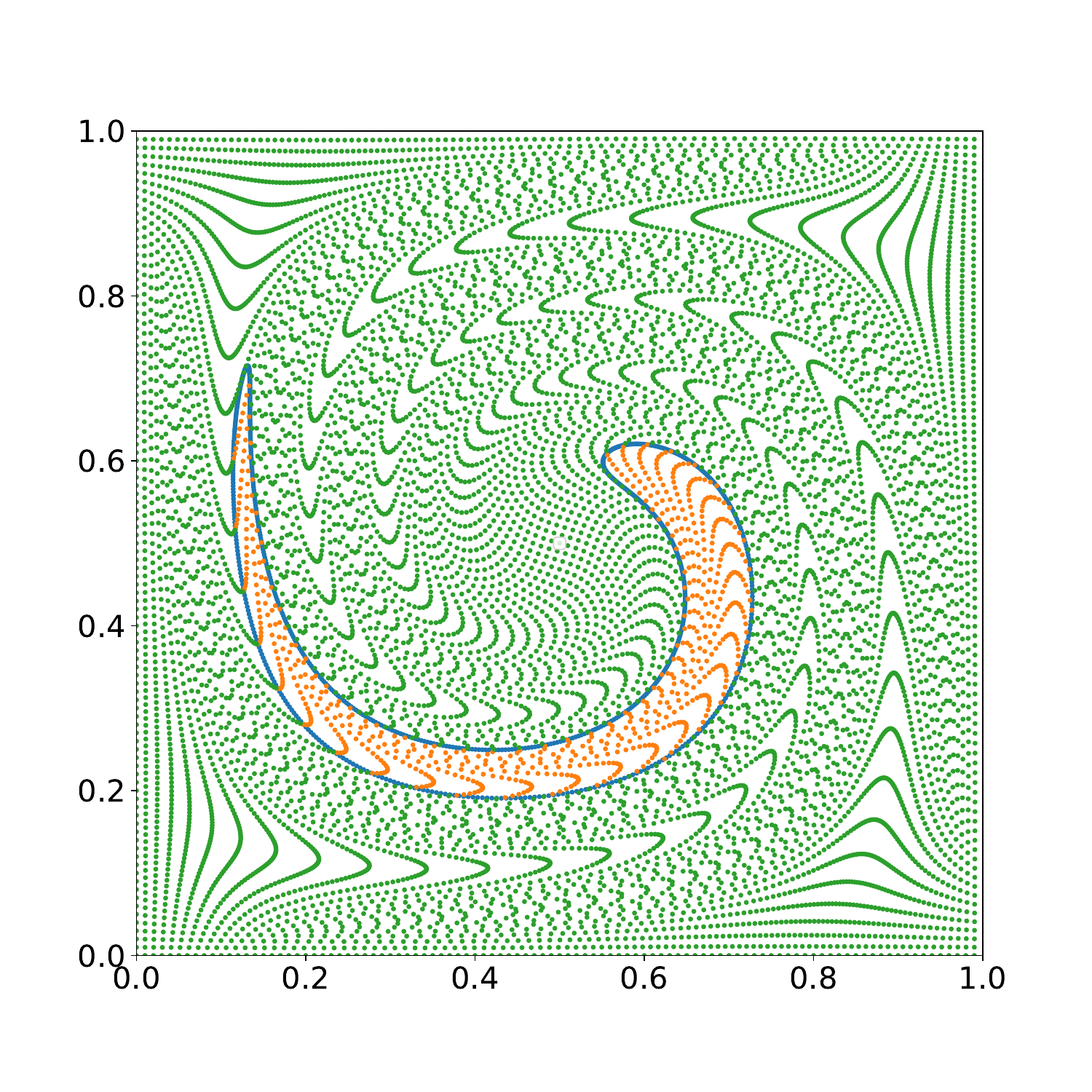}
	}
	\caption{Evolution of collocation points at $t = 0.0$, $0.5$, $1.0$, and $1.5$, respectively.}
	\label{fig3-2-2}
\end{figure}

In this problem, the flow velocity $\boldsymbol{w}$ is explicitly included in the equation, resulting in a large deformation.
For the implementation, we adopt the interface-tracking algorithm from \cite{doi:10.1137/17M1149328}. We utilize a $5^{th}$-order Runge-Kutta scheme for flow mapping approximation and the cubic MARS algorithm to capture large deformations in the domain. At the final time $T$, $\Omega_1(T)$ is stretched into a snake-like shape. Figure \ref{fig3-2-2} illustrates the evolution of collocation points at $t = 0.0$, $0.5$, $1.0$, and $1.5$, respectively.

The piecewise viscosity coefficients are set as $\nu_1=1$ and $\nu_2=10^{-3}$. The exact solution is determined by smooth velocity and pressure in each subdomain
\begin{equation*}
	\begin{aligned}
		u(\boldsymbol{x},t)=&\begin{cases}\cos t \cos \left(\pi x\right) \sin \left(\pi y\right), & \boldsymbol{x}\in\Omega_{1}(t) \\
			e^{x}\sin \left(\pi y+\pi t\right), & \boldsymbol{x}\in\Omega_{2}(t)\end{cases}, \\
		v(\boldsymbol{x},t)=&\begin{cases}-\cos t \sin \left(\pi x\right) \cos \left(\pi y\right), & \boldsymbol{x}\in\Omega_{1}(t) \\
			e^{x} \pi^{-1} \cos \left(\pi y+\pi t\right), & \boldsymbol{x}\in\Omega_{2}(t)\end{cases},\\			
		p(\boldsymbol{x},t)=&\begin{cases}\cos \left(0.5 \pi x\right) \sin \left(0.5 \pi y\right), & \boldsymbol{x}\in\Omega_{1}(t) \\
			\sin \left(0.5 \pi x\right) \cos \left(0.5 \pi y\right), & \boldsymbol{x}\in\Omega_{2}(t)\end{cases}.\\	
	\end{aligned}
\end{equation*}	

Table \ref{table3-2-2} records results of RFM and the high-order unfitted FEM proposed in \cite{ma2021high} for the Oseen equation. Due to the lack of discussion on the degrees of freedom in \cite{ma2021high}, the results for unfitted FEM in Table \ref{table3-2-2} only provide the mesh scale parameter $h$. However, it should be noted that in unfitted FEM, rectangular elements intersecting with the interface belong to both covers of two subdomains, resulting in a higher degree of freedom compared to RFM with the same $h$. Again, RFM achieves high accuracy even when solving simulations with significant deformations.

\begin{table}[htbp]
	\caption{\label{table3-2-2} Comparison of RFM and unfitted FEM for the Oseen equation.} \centering
	\begin{small}
		\begin{tabular}{|c|c|c|c|cc|}
			\hline
			Method & $h$ & $M$ & $N$ & $u$ error & $v$ error \\
			\hline
			\multirow{4}{*}{RFM}
			& \multirow{2}{*}{1/16} 
			& 28800 & 39916    &  1.01E-2   & 1.39E-2  \\
			& & 57600 & 39916  &  6.41E-3   & 6.17E-3  \\
			\cline{2-6}
			& \multirow{2}{*}{1/32} 
			& 57600 & 241020   &  4.51E-4   & 3.48E-4  \\
			& & 115200& 241020 &  9.34E-5   & 1.27E-4  \\
			\hline
			\multirow{2}{*}{unfitted FEM}
			& 1/16 & \multicolumn{2}{|c|}{——} &  \multicolumn{2}{|c|}{1.46E-3}  \\
			\cline{2-6}
			& 1/32 & \multicolumn{2}{|c|}{——} &  \multicolumn{2}{|c|}{8.70E-5}  \\
			\hline
		\end{tabular}
	\end{small}
\end{table}

\subsubsection{Linear fluid-solid interaction problem with complex geometry}
\label{sec3-2-3}
In order to account for the scenario where physical states differ between two subdomains, we consider a linear fluid-solid interaction problem introduced in \cite{kaltenbacher2018mathematical}, which can be formulated as follows
\begin{tiny}
\begin{equation}
	\begin{cases}
		\begin{aligned}
			\boldsymbol{u}^S_{t t}(\boldsymbol{x},t)-\Delta \boldsymbol{u}^S(\boldsymbol{x},t)+\boldsymbol{u}^S(\boldsymbol{x},t)+ \boldsymbol{u}^S_t(\boldsymbol{x},t) & \equiv \boldsymbol{f}^S(\boldsymbol{x},t), & & \quad  \boldsymbol{x} \in  \Omega_{1}(t), \quad t \in [0, T],\\
			\boldsymbol{u}^F_t(\boldsymbol{x},t)-\Delta \boldsymbol{u}^F(\boldsymbol{x},t)+\nabla p^F(\boldsymbol{x},t) & \equiv \boldsymbol{f}^F(\boldsymbol{x},t), && \quad  \boldsymbol{x} \in  \Omega_{2}(t), \quad t \in [0, T], \\
			\operatorname{div} \boldsymbol{u}^F(\boldsymbol{x},t) & \equiv 0, && \quad  \boldsymbol{x} \in  \Omega_{2}(t), \quad t \in [0, T], \\
			\boldsymbol{u}^F(\boldsymbol{x},t) & \equiv \boldsymbol{u}^S_t(\boldsymbol{x},t) ,&& \quad \boldsymbol{x} \in \Gamma(t), \quad t \in [0, T],\\
			\partial_{\boldsymbol{n}} \boldsymbol{u}^F(\boldsymbol{x},t)-\partial_{\boldsymbol{n}} \boldsymbol{u}^S_t(\boldsymbol{x},t) & =p^F(\boldsymbol{x},t) \boldsymbol{n}(\boldsymbol{x},t) ,&& \quad \boldsymbol{x} \in \Gamma(t), \quad t \in [0, T],\\
			\boldsymbol{u}^S(\boldsymbol{x},0)&=\boldsymbol{u}^S_{0}(\boldsymbol{x}), && \quad  \boldsymbol{x} \in \Omega_{1}(0) ,\\
			\boldsymbol{u}^S_t(\boldsymbol{x},0)&=\boldsymbol{u}^S_{1}(\boldsymbol{x}), && \quad  \boldsymbol{x} \in \Omega_{1}(0) ,\\
			\boldsymbol{u}^F(\boldsymbol{x},0)&=\boldsymbol{u}^F_{0}(\boldsymbol{x}), && \quad  \boldsymbol{x} \in \Omega_{2}(0) ,\\
			\boldsymbol{u}^F(\boldsymbol{x},t)&=\boldsymbol{g}(\boldsymbol{x}, t) ,&& \quad \boldsymbol{x} \in \partial \Omega, \quad t \in [0, T],
		\end{aligned}
	\end{cases}
	\label{pde-fsi}
\end{equation}
\end{tiny}
where $\boldsymbol{u}^S, \boldsymbol{u}^F$ denote the solid displacement and the flow velocity, respectively.

We consider problem \eqref{pde-fsi} over a fixed but complex interface geometry in Figure \ref{figure3-1-1-2} with the final time $T=1.0$.
And we specify the exact solid displacement and the flow velocity as
\begin{equation*}
	\begin{aligned}
			\boldsymbol{u}^S(\boldsymbol{x},t)&=\left(10-x^2-y^2-t^2, 20-x-y-t\right)^T ,&& \quad  \boldsymbol{x} \in \Omega_{1}(t)=\Omega_{1}, \\
			\boldsymbol{u}^F(\boldsymbol{x},t)&=\left(e^t \frac{y}{r},-e^t \frac{x}{r}\right)^T ,&& \quad  \boldsymbol{x} \in \Omega_{2}(t)=\Omega_{2},\\
			p^F(\boldsymbol{x},t)&=(x-1)^3+(y-1)^3+(t-1)^2 , && \quad  \boldsymbol{x} \in \Omega_{2}(t)=\Omega_{2}.
	\end{aligned}
\end{equation*}	

\begin{table}[htbp]
	\caption{\label{table3-2-3} Results of RFM for the fluid-solid interaction problem with complex interface geometry.} \centering
	\begin{small}
		\begin{tabular}{|c|c|cccc|}
			\hline
			$M$ & $N$ & $u^S$ error & $v^S$ error & $u^F$ error & $v^F$ error  \\
			\hline
			8000   & 58132  &  6.03E-5   & 7.19E-5  & 9.24E-4 & 1.29E-3 \\
			16000  & 58132  &  8.18E-6   & 2.04E-6  & 4.15E-5 & 3.15E-5 \\
			32000  & 58132  &  4.73E-6   & 2.31E-6  & 2.87E-6 & 2.39E-6 \\
			48000  & 58132  &  7.68E-6   & 8.20E-7  & 3.60E-7 & 6.99E-7 \\
			48000  & 323554 &  1.33E-9   & 4.61E-10 & 2.47E-8 & 2.15E-8 \\
			\hline
		\end{tabular}
	\end{small}
\end{table}
From Table \ref{table3-2-3}, it is evident that RFM still has high accuracy for this challenging problem.

\section{Conclusions}
\label{sec4}

In conclusion, our novel method effectively tackles the challenge of mesh generation for interface problems, which traditional methods struggle to overcome. This method is inspired by the observation that non-smooth solutions, which are also difficult to handle in interface problems, often exhibit piecewise smooth behavior. Therefore, we employ two sets of random feature functions to approximate the solution.  The interface conditions are seamlessly integrated into the loss function within the random feature method framework. Our method has been rigorously tested on a range of stationary and time-evolving interface problems of varying complexity. The results demonstrate that our method not only achieves high accuracy across all cases, but also significantly reduces the degrees of freedom by two to three orders of magnitude compared to traditional methods, while maintaining the same level of accuracy. Furthermore, our method negates the need for mesh generation, proving its robustness in handling time-dependent problems with complex interface geometries or intricate interface evolutions.

\section*{Acknowledgments}
The work is supported by National Key R\&D Program of China (No. 2022YFA1005200, No. 2022YFA1005202, and No. 2022YFA1005203), NSFC Major Research Plan -  Interpretable and General-purpose Next-generation Artificial Intelligence (No. 92270001 and No. 92270205), Anhui Center for Applied Mathematics, and the Major Project of Science \& Technology of Anhui Province (No. 202203a05020050). We thank Professor Weinan E for helpful discussions.

\begin{appendix}

\section{Experimental Setup}\label{sec::appendix}

Table \ref{table3-1} records the hyper-parameters in RFM for stationary interface problems in Section \ref{sec3-1}.
\begin{table}[htbp]
	\caption{\label{table3-1} Hyper-parameters in RFM for stationary interface problems.} \centering
	\begin{small}
		\begin{tabular}{|c|c|c|c|c|}
			\hline
			Interface problem & $M_p$ & \multicolumn{2}{|c|}{$J_n$} & $M$ \\
			\hline
			2D Elliptic (Table \ref{table3-1-1-1}) & 4 & \multicolumn{2}{|c|}{$200$, $400$}  &  $2 M_{p} J_{n}$\\
			2D Elliptic (Table \ref{table3-1-1-2}) & 64 & \multicolumn{2}{|c|}{$400$}  &  $2 M_{p} J_{n}$\\
			2D Stokes (Table \ref{table3-1-2-1}) & 16 & \multicolumn{2}{|c|}{$100$, $400$}  &  $6 M_{p} J_{n}$\\
			3D Stokes (Table \ref{table3-1-2-2}) & 8 & \multicolumn{2}{|c|}{$600$}  &  $8 M_{p} J_{n}$\\
			3D Stokes (Table \ref{table3-1-2-3}) & 8 & \multicolumn{2}{|c|}{$800$, $1600$, $3200$}  &  $4 M_{p} J_{n}$\\
			3D Elasticity (Table \ref{table3-1-3}) & 8 & \multicolumn{2}{|c|}{$800$, $1200$}  &  $6 M_{p} J_{n}$\\
			
			\hline
			Interface problem & \multicolumn{3}{|c|}{$\boldsymbol{x}_{n}$} & $\boldsymbol{r}_{n}$ \\
			\hline
			2D Elliptic (Table \ref{table3-1-1-1}) & \multicolumn{3}{|c|}{ $\left\{(-1,-1), (-1,1), (1,-1), (1,1)\right\}$} & $(1,1)$\\
			2D Elliptic (Table \ref{table3-1-1-2}) & \multicolumn{3}{|c|}{ $\left\{(\frac{i}{4},\frac{j}{4}), i=7,9,j=5,7\right\}$} & $(\frac{1}{4},\frac{1}{4})$\\
			2D Stokes (Table \ref{table3-1-2-1}) & \multicolumn{3}{|c|}{ $\left\{ (\frac{2i-1}{2}-2,\frac{2j-1}{2}-2), i, j = 1,\cdots,4\right\}$} & $(\frac{1}{2},\frac{1}{2})$\\
			3D Stokes (Table \ref{table3-1-2-2}) & \multicolumn{3}{|c|}{ $\left\{(\pm\frac{1}{2},\pm\frac{1}{2},\pm\frac{1}{2})\right\}$} & $(\frac{1}{2},\frac{1}{2},\frac{1}{2})$\\
			3D Stokes (Table \ref{table3-1-2-3}) & \multicolumn{3}{|c|}{ $\left\{(\pm\frac{3}{5},\pm\frac{3}{5},\pm\frac{3}{5})\right\}$} & $(\frac{3}{5},\frac{3}{5},\frac{3}{5})$\\
			3D Elasticity (Table \ref{table3-1-3}) & \multicolumn{3}{|c|}{ $\left\{ (\frac{i}{4},\frac{j}{4},\frac{k}{4}),i,j,k=1,3\right\}$} & $(\frac{1}{4},\frac{1}{4},\frac{1}{4})$\\

			\hline
			Interface problem	& \multicolumn{2}{|c|}{$Q$} & \multicolumn{2}{|c|}{$N$} \\
			\hline
			2D Elliptic (Table \ref{table3-1-1-1}) & \multicolumn{2}{|c|}{$1600$, $6400$, $14400$} & \multicolumn{2}{|c|}{$3084$, $9370$, $18854$}\\
			2D Elliptic (Table \ref{table3-1-1-2}) & \multicolumn{2}{|c|}{$160000$, $230400$, $313600$} & \multicolumn{2}{|c|}{$185496$, $259920$, $348038$}\\
			2D Stokes (Table \ref{table3-1-2-1}) & \multicolumn{2}{|c|}{$1600$, $25600$} & \multicolumn{2}{|c|}{$6801$, $84801$}\\
			3D Stokes (Table \ref{table3-1-2-2}) & \multicolumn{2}{|c|}{$64000$} & \multicolumn{2}{|c|}{$129088$}\\
			3D Stokes (Table \ref{table3-1-2-3}) & \multicolumn{2}{|c|}{$64000$} & \multicolumn{2}{|c|}{$129088$}\\
			3D Elasticity (Table \ref{table3-1-3}) & \multicolumn{2}{|c|}{$64000$} & \multicolumn{2}{|c|}{$344820$}\\

			\hline
		\end{tabular}
	\end{small}
\end{table}

Hyper-parameters in RFM for time-dependent problems in Section \ref{sec3-2} are shown in Table \ref{table3-2}.
\begin{table}[htbp]
	\caption{\label{table3-2} Hyper-parameters in RFM for time-dependent interface problems.} \centering
	\begin{small}
		\begin{tabular}{|c|c|c|c|c|}
			\hline
			Interface problem & $M_p$ & \multicolumn{2}{|c|}{$J_n$} & $M$ \\
			\hline
			Parabolic moving (Table \ref{table3-2-1-1}) & 8 & \multicolumn{2}{|c|}{$800$, $1600$}  &  $2 M_{p} J_{n}$\\
			Parabolic moving (Table \ref{table3-2-1-2}) & 8 & \multicolumn{2}{|c|}{$50$, $100$, $200$, $400$, $800$}  &  $2 M_{p} J_{n}$\\
			Osean dynamic (Table \ref{table3-2-2}) & 8 & \multicolumn{2}{|c|}{$600$, $1200$, $2400$}  &  $6 M_{p} J_{n}$\\
			Fluid-solid interaction (Table \ref{table3-2-3}) & 8 & \multicolumn{2}{|c|}{$200$, $400$, $800$, $1200$}  &  $5 M_{p} J_{n}$\\
			
			\hline
			Interface problem & \multicolumn{3}{|c|}{$\boldsymbol{x}_{n}$} & $\boldsymbol{r}_{n}$ \\
			\hline
			Parabolic moving (Table \ref{table3-2-1-1}) & \multicolumn{3}{|c|}{ $\left\{(\pm\frac{1}{2},\pm\frac{1}{2},\frac{1}{4}),(\pm\frac{1}{2},\pm\frac{1}{2},\frac{3}{4})\right\}$} & $(\frac{1}{2},\frac{1}{2},\frac{1}{4})$\\
			Parabolic moving (Table \ref{table3-2-1-2}) & \multicolumn{3}{|c|}{ $\left\{(\pm\frac{1}{2},\pm\frac{1}{2},\frac{1}{4}),(\pm\frac{1}{2},\pm\frac{1}{2},\frac{3}{4})\right\}$} & $(\frac{1}{2},\frac{1}{2},\frac{1}{4})$\\
			Osean dynamic (Table \ref{table3-2-2}) & \multicolumn{3}{|c|}{ $\left\{(\frac{i}{4},\frac{j}{4},\frac{k}{8}), i,j,=1,3,k=3,9\right\}$} & $(\frac{1}{4},\frac{1}{4},\frac{3}{8})$\\
			Fluid-solid interaction (Table \ref{table3-2-3}) & \multicolumn{3}{|c|}{ $\left\{(\frac{i}{4},\frac{j}{4},\frac{k}{4}), i=7,9,j=5,7,k=1,3\right\}$} & $(\frac{1}{4},\frac{1}{4},\frac{1}{4})$\\

			\hline
			Interface problem	& \multicolumn{2}{|c|}{$Q$} & \multicolumn{2}{|c|}{$N$} \\
			\hline
			Parabolic moving (Table \ref{table3-2-1-1}) & \multicolumn{2}{|c|}{$64000$, $216000$} & \multicolumn{2}{|c|}{$106000$, $375360$}\\
			Parabolic moving (Table \ref{table3-2-1-2}) & \multicolumn{2}{|c|}{$64000$} & \multicolumn{2}{|c|}{$141261$}\\
			Osean dynamic (Table \ref{table3-2-2}) & \multicolumn{2}{|c|}{$64000$, $216000$} & \multicolumn{2}{|c|}{$39916$, $241020$}\\
			Fluid-solid interaction (Table \ref{table3-2-3}) & \multicolumn{2}{|c|}{$8000$, $64000$} & \multicolumn{2}{|c|}{$58132$, $323554$}\\
			
			\hline
		\end{tabular}
	\end{small}
\end{table}

\end{appendix}
\bibliographystyle{siam}
\bibliography{references}
	
\end{document}